\documentclass[epsfig,latexsym,amsfonts,twoside]{article}
\usepackage{amssymb}
\usepackage{amsmath,amscd}

\def\part#1{\frac{\partial\phantom{#1}}{\partial#1}}
\newtheorem{thm}{Theorem}
\newtheorem{theorem}[thm]{Theorem}
\newtheorem{proposition}[thm]{Proposition}
\newtheorem{lemma}[thm]{Lemma}
\newtheorem{corollary}[thm]{Corollary}

\newenvironment{proof}{\begin{trivlist}\item[]{\bf Proof} }%
{\hfill $\Box$ \end{trivlist}}
\newenvironment{definition}{\begin{trivlist}\item[]{\bf Definition}\em }%
{\end{trivlist}}
\newenvironment{remark}{\begin{trivlist}\item[]{\bf Remark} }%
{\end{trivlist}}
\newenvironment{example}{\begin{trivlist}\item[]{\bf Example} }%
{\end{trivlist}}

\def\Z{\ifmmode{{\mathbb Z}}\else{${\mathbb Z}$}\fi}
\def\Q{\ifmmode{{\mathbb Q}}\else{${\mathbb Q}$}\fi}
\def\C{\ifmmode{{\mathbb C}}\else{${\mathbb C}$}\fi} 
\def\P{\ifmmode{{\mathbb P}}\else{${\mathbb P}$}\fi} 

\def\H{\ifmmode{{\mathrm H}}\else{${\mathrm H}$}\fi} 

\def\B{\ifmmode{{\mathcal B}}\else{${\mathcal B}$}\fi} 
\def\E{\ifmmode{{\mathcal E}}\else{${\mathcal E}$}\fi} 
\def\F{\ifmmode{{\mathcal F}}\else{${\mathcal F}$}\fi} 
\def\K{\ifmmode{{\mathcal K}}\else{${\mathcal K}$}\fi} 
\def\L{\ifmmode{{\mathcal L}}\else{${\mathcal L}$}\fi} 
\def\M{\ifmmode{{\mathcal M}}\else{${\mathcal M}$}\fi} 
\def\N{\ifmmode{{\mathcal N}}\else{${\mathcal N}$}\fi} 
\def\O{\ifmmode{{\mathcal O}}\else{${\mathcal O}$}\fi} 
\def\U{\ifmmode{{\mathcal U}}\else{${\mathcal U}$}\fi}
\def\V{\ifmmode{{\mathcal V}}\else{${\mathcal V}$}\fi}
\def\X{\ifmmode{{\mathcal X}}\else{${\mathcal X}$}\fi} 

\def\Br{\ifmmode{{\mathrm{Br}}}\else{${\mathrm{Br}}$}\fi} 
\def\OG{\ifmmode{\widetilde{\cal M}_4}\else{$\widetilde{\cal M}_4$}\fi} 
\def\D{\ifmmode{{\mathcal D}_{\mathrm{coh}}^b}\else{${{\mathcal
    D}_{\mathrm{coh}}^b}$}\fi}
\def\Shah{\ifmmode{\amalg\hspace*{-3.5pt}\amalg}\else{$\amalg\hspace*{-3.5pt}\amalg$}\fi}

\begin{document}

\title{Twisted Fourier-Mukai transforms for holomorphic symplectic
  fourfolds}
\author{Justin Sawon}
\date{February, 2005}
\maketitle

\begin{abstract}
We apply the methods of C{\u a}ld{\u a}raru to construct a twisted
Fourier-Mukai transform between a pair of holomorphic symplectic
four-folds. More precisely, we obtain an equivalence between the
derived category of coherent sheaves on a certain four-fold and the
derived category of twisted sheaves on its `mirror' partner. As
corollaries, we show that the two spaces are connected by a
one-parameter family of deformations through Lagrangian fibrations,
and we extend the original Fourier-Mukai transform to degenerations of
abelian surfaces.
\end{abstract}

\maketitle

\section{Introduction}

Matsushita~\cite{matsushita99} proved that a projective holomorphic
symplectic manifold can only be fibred by (holomorphic) Lagrangian
abelian varieties; his results in~\cite{matsushita99,matsushita00}, 
together with the results of Cho, Miyaoka, and 
Shepherd-Barron~\cite{cmsb02}, also strongly suggest that the base 
of the fibration must be projective space. In~\cite{sawon03} the author
reviewed what is known about such fibrations, and speculated on what
may be true. In particular, we hope to obtain a classification (up to 
deformation) of holomorphic symplectic manifolds via this approach.

A central problem is: {\em can a Lagrangian fibration which does not
  admit a global section be deformed to one that does?\/} This is the
motivation behind the present paper, which answers the question
affirmatively for a particular example
(Theorem~\ref{Lagrangian_deformation}) while introducing ideas which
should have wider applications.

We investigate holomorphic symplectic four-folds fibred by abelian
surfaces. Following Altman and Kleiman~\cite{ak80}, we associate to
such a fibration $X\rightarrow B$ its compactified relative Picard
scheme $P:=\overline{\mathrm{Pic}}^0(X/B)$, which is fibred over $B$
and admits a section. We regard $P$ as the dual fibration, and the
double-dual fibration is $X^0:=\overline{\mathrm{Pic}}^0(P/B)$. If the
singular fibres of $X$ are not too bad, both $P$ and $X^0$ are
well-defined and smooth. It can also happen that $X$ and $X^0$ are
locally isomorphic as fibrations, and then $X$ is a torsor over
$X^0$. Our goal is to study the relation between $X$ and $X^0$, and to
show that $X$ can be deformed to $X^0$ via Lagrangian fibrations.

The analogous situation for elliptic surfaces was studied by
Kodaira~\cite{kodaira63,kodaira64};  for higher dimensional elliptic
fibred varieties it is known as Ogg-Shafarevich theory
(see~\cite{dg94} for instance). In those cases we use the compactified
relative Jacobian $J:=\overline{\mathrm{Jac}}(X/B)$ of an elliptic
fibration $X\rightarrow B$ (since elliptic curves are self-dual, we
don't need to take the double-dual fibration). The group classifying
all torsors $X$ over $J$ is the (cohomological) analytic Brauer group
$\H^2(J,\O^*)$ of $J$. Recently C{\u a}ld{\u a}raru~\cite{caldararu00}
gave a more conceptual explanation for the appearance of the Brauer
group in this context: $J$ can be interpreted as a moduli space of
sheaves on $X$, and there is a holomorphic gerbe
$\beta\in\H^2(J,\O^*)$ obstructing the existence of a universal sheaf
for the moduli problem. Then $X$ may be regarded as $J$ `twisted' by
$\beta$. C{\u a}ld{\u a}raru also defined a twisted Fourier-Mukai
transform, an equivalence between the derived category of coherent
sheaves on $X$ and the derived category of sheaves twisted by
$\beta^{-1}$ on $J$. The existence of such an equivalence is
indicative of the close geometric relation between the spaces $X$ and
$J$.

Among the examples considered by C{\u a}ld{\u a}raru are elliptic K3
surfaces and Calabi-Yau three-folds. In this paper we construct a 
twisted Fourier-Mukai transform for holomorphic symplectic four-folds 
fibred by Lagrangian abelian surfaces
(Theorem~\ref{twisted_equivalence}). The compactified relative Picard 
scheme $P$ will be, a priori, a moduli space of sheaves on the 
four-fold $X$. We will show that it is smooth and holomorphic 
symplectic, by identifying it with another well-known holomorphic
symplectic four-fold. The double-dual fibration $X^0$ will also be a 
smooth holomorphic symplectic four-fold. Moreover, $X$ will be a torsor 
over $X^0$, and will therefore correspond to a gerbe in $\H^2(P,\O^*)$. 
By analyzing the space $\H^2(P,\O^*)$ of gerbes on $P$, we will show
that there is a one-parameter family of Lagrangian fibrations 
connecting $X$, which does not admit a global section, to $X^0$, which 
does (Theorem~\ref{Lagrangian_deformation}). It should be stressed
that our main interest is in obtaining geometric results such as this,
and the application of the more abstract theory is just a tool to help
understand the structure of the fibrations. 

The twisted Fourier-Mukai transform will be an equivalence between the 
derived category of $X$ and the twisted derived category of $P$. It can
be regarded as a family version of the Fourier-Mukai transform between
an abelian surface and its dual. The twist arises when one tries to 
`assemble' the fibrewise transforms into a global transform. As a
corollary, we extend the Fourier-Mukai transform of abelian surfaces to 
degenerations of abelian surfaces (Corollary~\ref{singular_Mukai}).

On hyperk{\"a}hler manifolds, a holomorphic Lagrangian fibration becomes 
a special Lagrangian fibration after a rotation of complex structures. 
We expect that our equivalence is a manifestation of homological mirror 
symmetry in the presence of a B-field, as applied to SYZ 
fibrations~\cite{syz96}, although the precise relation between 
Fourier-Mukai transforms and homological mirror symmetry is not yet
understood. Similar examples of mirror pairs of hyperk{\"a}hler 
manifolds involving B-fields were constructed by Hausel and 
Thaddeus~\cite{ht01}, though their examples are non-compact. We expect
that our methods will produce more compact examples in higher 
dimensions, and that these will be related to Hausel and Thaddeus'
examples via the construction of Donagi, Ein, and
Lazarsfeld~\cite{del97} (namely, the compactified Hitchin system is a
degeneration of the Beauville-Mukai system).

The paper is organized as follows. In Section 2 we review results of
Mukai, Bridgeland, Maciocia, and C{\u a}ld{\u a}raru on Fourier-Mukai
and twisted Fourier-Mukai transforms. In Section 3 we introduce a pair
of holomorphic symplectic four-folds which are fibred by abelian
varieties and collect together some results about them. In Section 4
we construct a twisted Fourier-Mukai transform relating the derived
category and twisted derived category of the pair of four-folds from
Section 3. This is followed by our applications.

This paper was begun during a visit to the Institut des Hautes {\'
E}tudes Scientifiques and finished at the University of Kyoto; the 
author is grateful for the hospitality he received at both those 
places. The author has benefited from conversations with many people 
on the topics presented here: he thanks them all, particularly 
Tom Bridgeland, Andrei C{\u a}ld{\u a}raru, Eduardo de Sequeira 
Esteves, Steven Kleiman, Manfred Lehn, Yoshinori Namikawa, and Richard
Thomas. The author was supported in part by NFS Grant \#0305865.

\section{FM and twisted FM transforms}

We begin by reviewing Mukai's work~\cite{mukai81} on integral
transforms between derived categories. We state Bridgeland's
criterion~\cite{bridgeland99} for when we get an equivalence of
categories, and the removable singularities result of Bridgeland and 
Maciocia~\cite{bm02}. Then we review twisted Fourier-Mukai transforms, 
as they appear in C{\u a}ld{\u a}raru's thesis~\cite{caldararu00}.

\subsection{Fourier-Mukai transforms}

Suppose $M$ is a moduli space of semi-stable sheaves on some given
space $X$. The moduli space is fine if there exists a universal
sheaf $\U$ on $X\times M$
$$\begin{array}{ccccc}
 & & {\U} & & \\
 & & \downarrow & & \\
 & & X\times M & & \\
\pi_X & \swarrow & & \searrow & \pi_M \\
X & & & & M. \\
\end{array}$$
Given a sheaf $\E$ on $M$, we can define a sheaf on $X$ by
pulling $\E$ back to $X\times M$, tensoring with $\U$, and pushing
down to $X$. This map extends to a functor, known as an {\em integral
transform\/}, between the bounded derived categories of coherent
sheaves on $M$ and $X$, which we denote 
$$\Phi^{\U}_{M\rightarrow X}:\D (M)\rightarrow\D (X)$$
$$\E^{\bullet}\mapsto {\bf R}\pi_{X*}(\U\stackrel{\bf
  L}{\otimes}\pi_M^*\E^{\bullet}).$$ 
Likewise, using the dual sheaf ${\U}^{\vee}$ we obtain a functor
$$\Phi^{{\U}^{\vee}}_{X\rightarrow M}:\D (X)\rightarrow\D (M).$$
These functors were investigated by Mukai~\cite{mukai81,mukai87} in
various cases: for example, when $X$ and $M$ are dual elliptic curves,
dual abelian varieties, or K3 surfaces. In certain situations we get
an equivalence of triangulated categories, and the functors are then
known as {\em Fourier-Mukai transforms\/}. This phenomenon is closely
related to dualities in physics, such as mirror symmetry, and may be
regarded as some kind of quantum symmetry between spaces.

The following criterion, developed by Mukai, Bondal-Orlov, and
Bridgeland, tells us precisely when an integral transform is a
Fourier-Mukai transform. Let ${\O}_m$ denote the skyscraper sheaf
supported at $m\in M$. Then $\U_m:=\Phi^{\U}_{M\rightarrow X}{\O}_m$ is
the sheaf on $X$ which the point $m\in M$ parametrizes.

\begin{theorem}[Bridgeland~\cite{bridgeland99}]
\label{bridgeland}
Suppose that $X$ and $M$ are smooth and of the same dimension. The
functor $\Phi^{\U}_{M\rightarrow X}$ is an equivalence of triangulated
categories if and only if
\begin{enumerate}
\item for all $m\in M$, $\U_m\otimes{\K}_X=\U_m$ and $\U_m$
is simple, i.e.
$$\mathrm{Hom}_X(\U_m,\U_m)=\C,$$
\item and for all integers $i$ and for all $m_1\neq m_2\in M$,
$$\mathrm{Ext}^i_X(\U_{m_1},\U_{m_2})=0.$$
\end{enumerate}
\end{theorem}

\begin{remark}
The conditions in the theorem are roughly akin to the requirement that
$\{\U_m\}_{m\in M}$ behave like an orthonormal basis with respect to
the $\mathrm{Ext}^{\bullet}$-pairing on $\D(X)$.
\end{remark}

\begin{example}
We can take $X$ to be a smooth elliptic curve $E$, and $M$ to be its 
dual $\hat{E}$ (the Jacobian of $E$), regarded as the moduli space of
degree zero line bundles on $E$. The Poincar{\'e} line bundle provides
a universal bundle. In higher dimensions, we can take an abelian
variety $A$, its dual $\mathrm{Pic}^0A$, and the Poincar{\'e} line
bundle. These are the original examples of Mukai~\cite{mukai81}.
\end{example}

If $X$ is a curve or surface, techniques have been developed to
determine whether a particular moduli space of sheaves $M$ on $X$ is
smooth. In higher dimensions, however, much less is known. A priori we
may have no way of knowing whether $M$ is smooth, which is why the
following result of Bridgeland and Maciocia is particularly useful.
More significantly for us, it weakens slightly the conditions we need 
to check in order to show that an integral transform is a Fourier-Mukai 
transform.

\begin{theorem}[Bridgeland-Maciocia~\cite{bm02}, Proposition 6.1]
\label{bridgeland_maciocia}
Suppose $X$ is a smooth projective variety of dimension $n$. Let $M$ 
be a fine moduli space of sheaves on $X$, with $M$ an irreducible
projective scheme of dimension $n$. Let $\U$ be a universal sheaf,
with $\U_m$ defined as before. Suppose that
\begin{enumerate}
\item for all $m\in M$, $\U_m\otimes{\K}_X=\U_m$ and $\U_m$
is simple,
\item for all $m_1\neq m_2\in M$,
$$\mathrm{Hom}_X(\U_{m_1},\U_{m_2})=0,$$
and the closed subscheme
$$\Gamma(\U):=\{(m_1,m_2)\in M\times
M|\mathrm{Ext}^i_X(\U_{m_1},\U_{m_2})\neq 0\mbox{ for some }i\in\Z\}$$
of $M\times M$ has dimension at most $n+1$.
\end{enumerate}
Then $M$ is smooth and
$$\Phi^{\U}_{M\rightarrow X}:\D(M)\rightarrow\D(X)$$
is an equivalence of categories. In particular, $\Gamma(\U)$ must be
the diagonal in $M\times M$.
\end{theorem}

\subsection{Twisted Fourier-Mukai transforms}

So far we have considered only fine moduli spaces. There are two
kinds of obstructions to the existence of a global universal sheaf:
firstly, universal sheaves may not exist locally, and secondly, the
local universal sheaves may not patch together into a global sheaf. We
will assume that all semi-stable sheaves are actually stable, in which
case the first of these obstructions can be avoided. The second
obstruction was studied by C{\u a}ld{\u a}raru~\cite{caldararu00}; the
main results of this subsection are quoted from his thesis.

Choose an open cover $\{M_i\}$ of $M$ such that there exists a local
universal sheaf $\U_i$ over $X\times M_i$ for all $i$. Since $M$
parametrizes stable sheaves, local universal sheaves always exist on
small enough open sets; indeed it suffices to take a Stein covering 
of $M$. (If we include strictly semi-stable sheaves this is no longer 
the case, as a point in the moduli space can represent a whole 
S-equivalence class of sheaves.)

Consider the restrictions of the sheaves $\U_i$ and $\U_j$ to the
overlap $X\times M_{ij}=X\times (M_i\cap M_j)$.
$$\begin{array}{ccccc}
{\U}_i & & & & {\U}_j \\
 & \searrow & & \swarrow & \\
 & & X\times M_{ij} & & \\
\pi_X & \swarrow & & \searrow & \pi_M \\
X & & & & M_{ij} \\
\end{array}$$
Let $m\in M_{ij}$ and let $\U_m$ be the sheaf on $X$ represented by
$m$. Since both $\U_i$ and $\U_j$ are universal sheaves over $X\times 
M_{ij}$, it follows that their restrictions to $X\times m$ are
isomorphic, as
$$\U_i|_{X\times m}\cong\U_m\cong\U_j|_{X\times m}.$$
Since $M_{ij}$ is a Stein open set, these combine to give an
isomorphism
$$\phi_{ij}:\U_i|_{X\times M_{ij}}\rightarrow\U_j|_{X\times M_{ij}}.$$ 
Since $\U_m$ is stable and hence simple, the isomorphism
$\U_i|_{X\times m}\cong\U_j|_{X\times m}$ is given by multiplication
by a non-zero number. Then an alternate way of formulating the above
statement is to say that there is a line bundle $\L_{ij}$ on $M_{ij}$
such that $\U_i|_{X\times M_{ij}}$ is equal (not just isomorphic) to 
$\pi_M^*\L_{ij}\otimes\U_j|_{X\times M_{ij}}$. (It is important that
we have equality rather than isomorphism, as $\L_{ij}$ is of course
isomorphic to the trivial bundle on $M_{ij}$.)

On a triple intersection $X\times M_{ijk}$ composition gives
$$\phi_{ki}\circ\phi_{jk}\circ\phi_{ij}:\U_i|_{X\times M_{ijk}}
\rightarrow\U_i|_{X\times M_{ijk}}.$$
On each $X\times m$ this is once again multiplication by a non-zero
number, which we denote $\beta_{ijk}(m)$. Thus the composition is
determined by a section $\beta_{ijk}\in\Gamma(M_{ijk},\O^*)$. It can
be shown that these sections give a 2-cocycle representing a
cohomology class $\beta\in\H^2(M,\O^*)$.

\begin{definition}
A (holomorphic) gerbe on $M$ up to isomorphism is an element of the
(cohomological) analytic Brauer group $\H^2(M,\O^*)$. 
\end{definition}

\begin{remark}
Up to isomorphism, a line bundle corresponds to an element of
$\H^1(M,\O^*)$, and hence a gerbe up to isomorphism may be regarded as
a higher dimensional analogue of a line bundle.
\end{remark}

Gerbes themselves can be defined in various ways (see
Hitchin~\cite{hitchin01} for instance). One description is simply as a
2-cocycle representing the element of $\H^2(M,\O^*)$. Another
description is as a collection of line bundles on two-fold
intersections $M_{ij}$, such as our $\L_{ij}$, satisfying various
properties. One of these properties is that
$\L_{ij}\otimes\L_{jk}\otimes\L_{ki}$ should be trivial, which in our
case follows from
$$\U_i|_{X\times M_{ijk}}=\pi_M^*(\L_{ij}\otimes\L_{jk}\otimes\L_{ki}) 
\otimes\U_i|_{X\times M_{ijk}}.$$
Choosing trivializations of each $\L_{ij}$ then gives another
trivialization of $\L_{ij}\otimes\L_{jk}\otimes\L_{ki}$ which is given
by the section $\beta_{ijk}\in\Gamma(M_{ijk},\O^*)$. 

Note also that if we began with different local universal bundles 
$$\U_i^{\prime}=\pi_M^*\M_i\otimes\U_i,$$
where $\M_i$ are line bundles on $M_i$, then we would have obtained
different line bundles 
$$\L_{ij}^{\prime}=\M_i\otimes\M_j^{-1}\otimes\L_{ij}.$$
In this case $\{\L_{ij}\}$ and $\{\L_{ij}^{\prime}\}$ define
isomorphic gerbes (this is the definition of isomorphism, and is the
same as saying that two 2-cocycles represent the same class in
$\H^2(M,\O^*)$). Now if some choice leads to $\L_{ij}=\O_{M_{ij}}$ for
all $i$ and $j$, then the local universal sheaves agree on overlaps
and can be patched together to give a global universal sheaf.

\begin{proposition}[C{\u a}ld{\u a}raru~\cite{caldararu00}] Suppose 
that $M$ be a moduli space of stable sheaves on $X$. There is a gerbe 
$\beta\in\H^2(M,\O^*)$ (defined up to isomorphism) representing the
obstruction to the existence of a global universal sheaf on $X\times
M$. This is the sole obstruction: $\beta$ vanishes if and only if
there exists a global universal sheaf.
\end{proposition}

\begin{remark}
Henceforth everything we will say can be made to depend only on the
isomorphism class of the gerbe, and therefore we will refer to gerbes
up to isomorphism simply as gerbes.
\end{remark}

Since the obstruction is completely encoded in the gerbe $\beta$, we
can construct Fourier-Mukai transforms by incorporating $\beta$ into
the construction. Specifically, this means working with {\em twisted
sheaves\/}.

\begin{definition}
Let $\beta$ be a gerbe on $M$. A $\beta$-twisted sheaf on $M$ is a
collection of sheaves $\F_i$ on $M_i$ and isomorphisms
$\psi_{ij}:\F_i|_{M_{ij}}\rightarrow\F_j|_{M_{ij}}$ on the overlaps
$M_{ij}$ such that
\begin{enumerate}
\item for all $i$ and $j$, $\psi_{ji}=\psi_{ij}^{-1}$,
\item and for all $i$, $j$, and $k$, the composition
$$\psi_{ki}\circ\psi_{jk}\circ\psi_{ij}:
\F_i|_{M_{ijk}}\rightarrow\F_{M_{ijk}}$$
is given by $\beta_{ijk}\mathrm{Id}$.
\end{enumerate}
\end{definition}

\begin{example}
The gerbe $\beta\in\H^2(M,\O^*)$ can be pulled back by $\pi_M$ to give
a gerbe $\pi_M^*\beta$ on $X\times M$. The collection $\{U_i\}$ of
local universal sheaves with isomorphisms $\phi_{ij}$ then gives a 
$\pi_M^*\beta$-twisted sheaf on $X\times M$. Let us denote this {\em
twisted universal sheaf\/} simply by $\U$, as in the untwisted case.
\end{example}

The category of $\beta$-twisted sheaves over $M$ is an abelian
category, and one can construct its derived category $\D(M,\beta)$. As
in the untwisted case, we can construct a functor
$$\Phi^{\U}_{M\rightarrow X}:\D (M,\beta^{-1})\rightarrow\D (X)$$
$$\E^{\bullet}\mapsto {\bf R}\pi_{X*}(\U\stackrel{\bf
  L}{\otimes}\pi_M^*\E^{\bullet}).$$
The inverse of $\beta$ is defined in the obvious way
$$(\beta^{-1})_{ijk}=(\beta_{ijk})^{-1}.$$
Note that $\E^{\bullet}$ is a complex of $\beta^{-1}$-twisted sheaves,
so $\pi_M^*\E^{\bullet}$ is a complex of $\pi_M^*\beta^{-1}$-twisted
sheaves; when we tensor it with $\U$, which is a
$\pi_M^*\beta$-twisted sheaf, the twistings cancel each other, and
hence
$$\U\stackrel{\bf L}{\otimes}\pi_M^*\E^{\bullet}$$ 
is a (untwisted) sheaf on $X\times M$. 

If $\Phi^{\U}_{M\rightarrow X}$ is an equivalence of triangulated
categories it is called a {\em twisted Fourier-Mukai transform\/}. C{\u
 a}ld{\u a}raru generalized Bridgeland's criterion for when this
happens. First observe that the skyscraper sheaf $\O_m$ on $M$ can be
regarded as a twisted sheaf; simply choose the cover $\{M_i\}$ so that
$m$ lies in precisely one open set, then just one local sheaf is
non-vanishing and all isomorphisms are zero. As in the untwisted case,
$\U_m:=\Phi^{\U}_{M\rightarrow X}\O_m$ is the sheaf parametrized by
the point $m\in M$.

\begin{theorem}[C{\u a}ld{\u a}raru~\cite{caldararu00}, Theorem 3.2.1]
\label{equivalence_twisted}
Suppose that $X$ and $M$ are smooth and of the same dimension, and
$\U$ is a $\pi_M^*\beta$-twisted universal sheaf on $X\times M$. The
functor $\Phi^{\U}_{M\rightarrow X}$ is an equivalence of triangulated
categories if and only if
\begin{enumerate}
\item for all $m\in M$, $\U_m\otimes{\K}_X=\U_m$ and $\U_m$
is simple, i.e.
$$\mathrm{Hom}_X(\U_m,\U_m)=\C,$$
\item and for all integers $i$ and for all $m_1\neq m_2\in M$
$$\mathrm{Ext}^i_X(\U_{m_1},\U_{m_2})=0.$$
\end{enumerate}
\end{theorem}

\subsection{Elliptic fibrations}

C{\u a}ld{\u a}raru~\cite{caldararu00} discusses two main examples of 
twisted Fourier-Mukai transforms: when $X$ and $M$ are K3 surfaces, 
and when $X$ and $M$ are elliptic Calabi-Yau three-folds. We will 
focus on the second case, and then generalize some of the results to 
fibrations by abelian varieties.

Let $p_X:X\rightarrow B$ be an elliptic fibration of arbitrary
dimension. Assume for the moment that all the fibres are reduced and
irreducible (i.e.\ either they are smooth elliptic curves, or they
contain a single node or cusp). 

\begin{definition}[D'Souza~\cite{dsouza79}]
The compactified relative Jacobian $J:=\overline{\mathrm{Jac}}(X/B)$
of $X$ is the moduli space parametrizing families of torsion-free
rank one sheaves of degree zero on the fibres of $X$. The degree of a 
sheaf $\E$ on a fibre $X_t$ is defined by $\chi(\E)-\chi(\O_{X_t})$.
\end{definition}

\begin{remark}
\label{identify}
Given a degree zero torsion-free rank one sheaf on a fibre, we can 
push-forward by the inclusion of the fibre in $X$, to obtain a torsion 
sheaf on $X$ itself. According to Simpson's 
terminology~\cite{simpson94}, this sheaf is pure-dimensional. Since 
it is rank one on a fibre of $X$, and the fibres are irreducible, a 
destabilizing sheaf cannot exist. A little more work shows that this
actually gives an isomorphism between $J$ and an irreducible component
of the Simpson moduli space of stable sheaves on $X$, and we will
henceforth use these two descriptions interchangeably.
\end{remark}

There is the obvious projection $p_J:J\rightarrow B$, and $J$ is
locally isomorphic to $X$ as a fibration: it is clear that
corresponding smooth fibres of $X$ and $J$ are isomorphic, and in fact
this is true also for singular fibres (see Section 6.3 of C{\u a}ld{\u
 a}raru~\cite{caldararu00}). Smoothness of $J$ now follows. Note that
if we allowed $X$ to have fibres with worse singularities then $J$ need 
not be smooth. The case of a type $I_2$ singular fibre is studied in 
Section 6 of C{\u a}ld{\u a}raru~\cite{caldararu00}; since it is not
irreducible, the elements of $J$ supported on such a fibre need not be
stable. This leads to singularities in $J$. 

Now $J$ can be regarded as the family of elliptic curves dual to the
family $X\rightarrow B$. The fact that $X$ and $J$ are locally
isomorphic is then a consequence of the self-duality of elliptic
curves (in higher dimensions it won't always be true, as abelian
varieties are not always self-dual). So we have a family version of an
elliptic curve and its dual. Locally we can give an explicit
construction of universal sheaves by extending the Poincar{\'e} line
bundle to a family over a small open set $B_i\subset B$. Note that
this actually gives a local universal bundle on the fibre product
$X\times_{B_i}J_i$, where $J_i:=p_J^{-1}(B_i)$, but it can be push
forward to $X\times J_i$ by the natural embedding. As usual, these
local universal sheaves need not patch together to give a global
universal sheaf: the obstruction is a gerbe $\beta\in\H^2(J,\O^*)$.

In this example $J$ admits a canonical section (given by the trivial
line bundle on each fibre of $X$) and $X$ is a torsor over
$J$. Thus $X$ is isomorphic to $J$ if and only if it admits a global
section. Moreover, to extend the Poincar{\'e} line bundle to a local
family is essentially the same as choosing a local section of $X$; if
we have a global section it can be extended globally. Such arguments 
lead to the next result.

\begin{proposition}[C{\u a}ld{\u a}raru~\cite{caldararu00}]
The following are equivalent
\begin{enumerate}
\item $X$ is isomorphic to $J$,
\item $X$ admits a section,
\item there is a global universal sheaf on $X\times J$,
\item $\beta$ vanishes.
\end{enumerate}
\end{proposition}

Moreover, when $\beta$ does not vanish, the fibration $X$ can be
reconstructed from $J$ and $\beta$ (this construction will be
described in the next subsection). If $J$ is projective, then $X$ is
also projective if and only if the corresponding gerbe $\beta$ is a
torsion element in $\H^2(M,\O^*)$, or equivalently, an element of the
{\'e}tale cohomology group $\H^2_{\mathrm{\acute{e}t}}(J,\O^*)$ known
as the (cohomological) {\em Brauer group\/}. Thus C{\u a}ld{\u
  a}raru's approach gives a conceptual explanation for the appearance
of the Brauer group in Ogg-Shafarevich theory~\cite{dg94}, where it is
essentially the group classifying all (minimal, projective) elliptic 
fibrations with a given relative Jacobian $J$.

There are of course some subtleties: in general elliptic fibrations
will have fibres which are not irreducible, nor reduced. For instance,
it is not currently known whether there exists an elliptic Calabi-Yau
three-fold with {\em only\/} irreducible fibres. Also, the following 
remark shows that not any element of the (cohomological) Brauer group 
can be used to construct a fibration $X$ from $J$.

\begin{remark}
Gerbes on $B$ can be pulled back to $J$ by the projection $p_J$,
giving an inclusion
$$\H^2(B,\O^*)\hookrightarrow\H^2(J,\O^*).$$
Moreover $J\rightarrow B$ admits a section so the inclusion splits. A
refinement of the theory shows that the obstruction $\beta$ really
lies in $\H^2(J,\O^*)/\H^2(B,\O^*)$.
\end{remark}

\subsection{Abelian fibrations}

The main purpose of this article is to describe twisted Fourier-Mukai
transforms for fibrations by abelian surfaces. We described an
elliptic fibration and its relative Jacobian as a family version of an
elliptic curve and its dual. We now want to discuss the family version
of a higher dimensional abelian variety and its dual.

Let $p_X:X\rightarrow B$ be a fibration by abelian varieties. Once
again, difficulties arise for fibres which are `too' singular, so we
will assume that all fibres are reduced and irreducible (this 
assumption will be satisfied by our main example later on).

\begin{definition}[Altman-Kleiman~\cite{ak80}]
The compactified relative Picard scheme
$P:=\overline{\mathrm{Pic}}^0(X/B)$ of $X$ is the moduli space
parametrizing families of torsion-free rank one sheaves of degree zero
on fibres of $X$. We say a sheaf $\E$ on a fibre $X_t$ has degree zero
if it is in the same connected component of 
$\overline{\mathrm{Pic}}(X_t)$ as $\O_{X_t}$.
\end{definition}

\begin{remark}
As with the compactified relative Jacobian of an elliptic fibration, 
we can take the push-forward of these sheaves under the inclusion of 
the fibre, and hence regard them as torsion sheaves on $X$ itself.
Then they are sheaves of pure-dimension according to 
Simpson~\cite{simpson94}. Since we have assumed the fibres of $X$ are
reduced and irreducible, there can be no destabilizing sheaves. Thus
$P$ parametrizes a family of stable sheaves on $X$. 

The typical element in this family is the push-forward of a degree
zero line bundle supported on a fibre of $X$. Let $Q$ be the
irreducible component of the Simpson moduli space of stable sheaves
on $X$ containing this typical element. One can show that we have a
(holomorphic) embedding $P\hookrightarrow Q$. The following lemma
implies that $P$ is actually isomorphic to $Q$, an important fact for
our later applications. Thus we can (and will) identify the
compactified relative Picard scheme $P$ with a component of the
Simpson moduli space of stable sheaves on $X$.
\end{remark}

Note that every sheaf in $Q$ is supported on a fibre, and thus we have
a projection $p_Q:Q\rightarrow B$.

\begin{lemma}
\label{no_fat_sheaves}
Every point $m$ of $Q$ parametrizes a sheaf $\U_m$ on $X$ of 
the form $\iota_*\E$, where $\iota:X_t\hookrightarrow X$ is the
inclusion of a fibre ($t:=p_Q(m)\in B$) and $\E$ is a stable rank one
degree zero sheaf on $X_t$.
\end{lemma}

\begin{proof}
This essentially follows from the assumption that the fibres of $X$
are reduced and irreducible. As above, this implies that all
sheaves in $Q$ are stable, and hence local universal sheaves
exist. Let $\U_i$ be a local universal sheaf over $X\times Q_i$, where 
$Q_i\subset Q$ is some open subset. Let $\Gamma$ be the graph of the
projection $p_X:X\rightarrow B$ in $X\times B$. We claim
that $\U_i$ is supported on the inverse image of $\Gamma$ under the
map given by the identity on $X$ and projection on $Q_i$.
$$\begin{array}{ccc}
{\U}_i & & \\
\downarrow & & \\
X\times Q_i & \stackrel{\mathrm{Id}_X\times
  p_Q}{\longrightarrow} & X\times B \\
\cup & & \cup \\
(\mathrm{Id}_X\times p_Q)^{-1}(\Gamma) & & \Gamma \\
\end{array}$$
The inverse image of $\Gamma$ and the support of $\U_i$ are closed
subsets cut out locally by some set of algebraic equations, and hence
it suffices to check that they agree near smooth fibres. But this is
obvious, as the sheaves $\U_m$ supported on a smooth fibre
$X_t$ (i.e.\ a smooth abelian variety) are just push-forwards of line
bundles on the fibre, and hence the claim follows.

Now if $Q_t$ is a fibre of $Q_i\subset P$ (where $t\in B$), then it
follows from the above claim that the restriction of $\U_i$ to
$X\times Q_t$ must be supported on $X_t\times Q_t$. Since
$X_t\subset X$ is irreducible and reduced, it follows that
$\U_i|_{X\times Q_t}$ is the push-forward of a sheaf on $X_t\times
Q_t$. The lemma now follows. 
\end{proof}

\begin{remark}
Although the proof relied on the existence of a local universal sheaf,
which does not exist when $X$ contains reducible fibres (and hence $Q$
contains strictly semi-stable sheaves), one expects the result to
remain true provided no fibres of $X$ contain non-reduced
components. On the other hand, if $X$ contained non-reduced fibres, 
we expect there could exist `fat' sheaves supported on those fibres,
i.e.\ sheaves which are not the push-forward of a sheaf on the
reduction of the fibre. The point of the lemma is to avoid such 
behaviour.
\end{remark}

Unlike in the case of an elliptic fibration, $P$ need not be locally
isomorphic to $X$ as a fibration, for two reasons. Firstly, if the
fibres of $X$ are not principally polarized then even a smooth fibre
of $P$ need not be isomorphic to the corresponding smooth fibre of $X$
(they will only be isogenous). Secondly, even for principal
polarizations the corresponding singular fibres of $P$ and $X$ may not
be isomorphic.

Regarding the second of these problems, let $X_t$ be a singular fibre
of $X$. We have a description of the singular fibre $P_t$ of $P$ as 
the moduli space of stable rank one degree zero sheaves on $X_t$. 
Further analysis needs to be done on a case-by-case basis, depending 
on the structure of $X_t$.

The first problem can be resolved by taking the double dual of $X$,
namely $X^0:=\overline{\mathrm{Pic}}^0(P/B)$. Clearly this is locally
isomorphic to $X$ as a fibration away from singular fibres. We will
postpone further discussion of singular fibres to specific examples; 
at this stage let us just say that the local isomorphisms can 
sometimes be extended over singular fibres, and we will assume this 
in the following discussion.

\begin{remark}
To construct $X^0$ from $X$ in one step, we can take the relative
Albanese scheme of $X$, following
Markushevich~\cite{markushevich96}. However, use of the relative
Picard scheme is essential for our moduli space interpretation.
\end{remark}

Take a (Stein) open cover $\{B_i\}$ of the base $B$, such that we 
have (local) isomorphisms
$$f_i:X_i\rightarrow X^0_i$$
where $X_i:=p_X^{-1}(B_i)$ and $X^0_i:=p_0^{-1}(B_i)$. On the overlap
$X^0_{ij}:=p_0^{-1}(B_{ij})$, define $\alpha_{ij}:=f_j\circ
f_i^{-1}$. Then 
$$\alpha_{ij}:X^0_{ij}\rightarrow X^0_{ij}$$
is given by a translation in each fibre over $B_{ij}$. Since $X^0$ has
a canonical section, a translation is equivalent to a local section of
$p^0:X^0\rightarrow B$. These local sections form a 1-cocycle
$\alpha\in\H^1(B,X^0)$, where we have implicitly identified $X^0$ with
its sheaf of local sections. Moreover $X$ can be completely recovered
from $\alpha$, so there is a bijection between torsors $X$ over $X^0$
and elements of $\H^1(B,X^0)$. 

By definition each point of $X^0$ represents a stable rank one
degree zero sheaf (generically a line bundle) supported on a fibre of
$P\rightarrow B$. Therefore each local section $\alpha_{ij}$ of
$p^0:X^0\rightarrow B$ determines a line bundle $\L_{ij}$ on
$P_{ij}:=p_P^{-1}(B_{ij})$, and the collection of these line bundles
determines the gerbe $\beta\in\H^2(P,\O^*)$. Thus we have shown the
next result.

\begin{proposition}
\label{abelian_prop}
The following are equivalent
\begin{enumerate}
\item $X$ is isomorphic to $X^0$,
\item $X$ admits a section,
\item there is a global universal sheaf on $X\times P$,
\item $\beta$ vanishes.
\end{enumerate}
\end{proposition}

Conversely, if the gerbe $\beta\in\H^2(P,\O^*)$ can be represented by
a collection of line bundles $\L_{ij}$ on $P_{ij}$ which have degree
zero when restricted to each fibre of $P$, then the construction can
be reversed. Thus in this case, there exists a torsor $X$ over $X^0$
corresponding to the gerbe $\beta$.

\section{Holomorphic symplectic manifolds}

We review some examples of holomorphic symplectic manifolds. In
particular, we describe some holomorphic symplectic four-folds which
are fibred by abelian surfaces, and collect together some facts about
these spaces that will be used in the next section.

\subsection{Definition and examples}

\begin{definition}
Let $X$ be a compact K{\"a}hler manifold. We call $X$ a holomorphic
symplectic manifold if it admits a closed non-degenerate two-form 
$\sigma$ of type $(2,0)$, i.e.\ 
$$\sigma\in\H^0(X,\Lambda^2T^*)=\H^{2,0}(X),$$
which we call a holomorphic symplectic form. If $X$ is
simply-connected and $\sigma$ generates $\H^{2,0}(X)\cong\C$ then we
say that $X$ is irreducible.
\end{definition}

If $X$ has (complex) dimension $2n$ then $\sigma^{\wedge n}$
trivializes the canonical bundle $\K_X$. By Yau's theorem, $X$ admits a
hyperk{\"a}hler metric; conversely, a hyperk{\"a}hler manifold is
holomorphic symplectic for each choice of complex structure compatible
with the hyperk{\"a}hler metric. By the Bogomolov decomposition theorem,
a holomorphic symplectic manifold has a finite cover which is the
cartesian product of a complex tori and irreducible holomorphic
symplectic manifolds. In this sense, all holomorphic symplectic
manifolds can be built out of irreducible ones.

In dimension two, K3 surfaces are the only irreducible examples, and
they form a single family up to deformation. In dimension four, there
are just two currently known examples, up to deformation.

\begin{example}
The first higher dimensional example was discovered by
Fujiki~\cite{fujiki83}. Let $S$
be a K3 surface and $\mathrm{Blow}_{\Delta}(S\times S)$ the blow up of
the diagonal. Quotienting by the involution which exchanges the two
copies of $S$ gives a smooth four-fold 
$$\mathrm{Hilb}^2S:=\mathrm{Blow}_{\Delta}(S\times S)/\Z_2.$$
Fujiki showed that $\mathrm{Hilb}^2S$ is an irreducible holomorphic
symplectic four-fold. Beauville~\cite{beauville83} generalized this
example to produce an irreducible holomorphic symplectic manifold in
each even dimension $2n$. These are the Hilbert schemes
$\mathrm{Hilb}^nS$, which parametrizes length $n$ zero-dimensional
subschemes of $S$, and are smooth resolutions of the symmetric
products $\mathrm{Sym}^nS$.
\end{example}

By beginning with an abelian surface, instead of a K3 surface, 
Beauville~\cite{beauville83} also constructed another family of
examples, one in each even dimension, known as the generalized Kummer
varieties. The following examples will also be important, though up to
deformation they do not give us new spaces.

\begin{example}
Let $S$ be a K3 surface with ample divisor $H$. The Mukai lattice is 
$$\H^{\bullet}(S,\Z)=\H^0(S,\Z)\oplus\H^2(S,\Z)\oplus\H^4(S,\Z)$$
endowed with the bilinear form
$$((v_0,v_2,v_4),(w_0,w_2,w_4)):=\int_S -v_0w_4+v_2w_2-v_4w_0.$$
The Mukai vector of a sheaf $\E$, defined by
$v(\E):=\mathrm{ch}(\E)\mathrm{td}^{1/2}\in\H^{\bullet}(S,\Z)$, is a
convenient way to encode the topological type of the sheaf. For
example, if $\E$ is a rank $r$ vector bundle with Chern classes $c_1$
and $c_2$ then
$$v(\E)=(r,c_1,r+c_1^2/2-c_2).$$
For fixed $v$ in the Mukai lattice, the Mukai moduli space $\M^s_H(v)$ 
is the moduli space of stable (with respect to $H$) sheaves $\E$ on
$S$ with fixed Mukai vector $v(\E)=v$.

Mukai~\cite{mukai84} showed that, for general $H$, $\M^s_H(v)$ is
smooth, quasi-projective, and holomorphic symplectic of dimension
$2n:=(v,v)+2$. If $v$ is primitive and $v_0>0$ then $\M^s_H(v)$ is
also compact: in fact it is an irreducible holomorphic symplectic
manifold. It is also a deformation of $\mathrm{Hilb}^nS$ (these 
results were proved by G{\"o}ttsche, Huybrechts, O'Grady, and
Yoshioka; see~\cite{yoshioka99}). For other choices of $v$, we can
compactify to the moduli space of semi-stable sheaves $\M^{ss}_H(v)$,
but this introduces singularities. 
\end{example}

\subsection{Abelian fibrations}

Elliptic K3 surfaces are dense and of codimension one in the moduli
space of all K3 surfaces. In higher dimensions we have the following
result.

\begin{theorem}[Matsushita~\cite{matsushita99,matsushita00}]
Let $X^{2n}$ be a projective irreducible holomorphic symplectic
manifold. Suppose $p_X:X\rightarrow B$ is a proper surjective
morphism, whose generic fibre is connected, and with projective base 
$B$ of dimension strictly between $0$ and $2n$. Then 
\begin{enumerate}
\item the generic fibre is a (holomorphic) Lagrangian abelian variety
  of dimension $n$,
\item the base is Fano with the same Hodge numbers as $\P^n$.
\end{enumerate}
In particular, when $n=2$ the base is $\P^2$.
\end{theorem}

It is currently an open problem whether an arbitrary holomorphic
symplectic manifold can be deformed to a fibration by abelian
varieties. This is possible for all the known irreducible 
holomorphic symplectic manifolds. In particular, for 
$\mathrm{Hilb}^nS$ we have the following example.

\begin{example}
Let $S$ be a K3 surface which contains a smooth genus $g\geq 2$ curve 
$C$, but is otherwise generic. Then $C$ moves in a $g$-dimensional
linear system $|C|\cong\P^g$, and taking $H=C$ as an ample divisor
gives us an embedding $S\hookrightarrow(\P^g)^{\vee}$ (unless $g=2$, 
in which case we instead get a double cover of the plane). Let $Z$ be
the Mukai moduli space $\M^s_H((0,[C],1))$, where $[C]$ denotes the
class of $C$ in $\H^2(S,\Z)$. Then $Z$ is smooth and compact. The
typical element is the push-forward of a degree $g$ line bundle on a
smooth curve $D\in |C|$. Thus $Z$ is fibred over $|C|\cong\P^g$, and
the generic fibre is a smooth abelian variety of dimension $g$, namely
the degree $g$ Picard group $\mathrm{Pic}^gD$ of a smooth genus $g$
curve. 
$$\begin{array}{ccc}
\mathrm{Pic}^g & \hookrightarrow & Z \\
               &                 & \downarrow \\
               &                 & |C|\cong\P^g \\
\end{array}$$

Let us argue that $Z$ is birational to $\mathrm{Hilb}^gS$. A generic
element of $Z$ gives a generic degree $g$ line bundle on a smooth
genus $g$ curve $D\subset S$. This line bundle will have a unique
section, up to scale, which will vanish at precisely $g$ distinct
points. This gives us a rational map $Z\dashrightarrow\mathrm{Hilb}^gS$.

In the other direction, a generic element of $\mathrm{Hilb}^gS$
consisting of $g$ distinct points determines a hyperplane in
$(\P^g)^{\vee}$. For $g>2$ this hyperplane cuts
$S\subset(\P^g)^{\vee}$ in a smooth curve of the linear system $|C|$;
for $g=2$ it can be pulled back from $(\P^2)^{\vee}$ to give such a
curve. Moreover, the $g$ points lie on the curve and determine a
degree $g$ line bundle. Thus we obtain a rational map
$\mathrm{Hilb}^gS\dashrightarrow Z$.
\end{example}

This example was used by Beauville~\cite{beauville99} to count the
number of rational curves (with nodes) in each linear system $|C|$.
Note that as in the first remark in Section~\ref{identify}, $Z$ can be
identified with the compactified relative Jacobian
$\overline{\mathrm{Pic}}^g(\mathcal{C}/|C|)$ (see
D'Souza~\cite{dsouza79} or Altman, Iarrobino, and
Kleiman~\cite{aik77}) of the family of curves $\mathcal{C}\rightarrow
|C|$ in the linear system. Huybrechts~\cite{huybrechts97} showed that
birational holomorphic symplectic manifolds have the same periods and
represent non-separated points in the moduli space, which is
non-Hausdorff. It follows that $Z$ and $\mathrm{Hilb}^gS$ are also
deformation equivalent.

Debarre~\cite{debarre99} used similar methods to Beauville to count
the number of genus two curves (with nodes) in a linear system on an
abelian surface. This included showing that the generalized Kummer
varieties can be deformed to fibrations by abelian varieties (see also
Example 3.8 in Sawon~\cite{sawon03}).

\subsection{More about $Z$}

In this subsection we concentrate on the fibrations by abelian
surfaces which are deformations of $\mathrm{Hilb}^2(S)$. We collect
together some facts that will be of use in the next section.  

Let $S\rightarrow(\P^2)^{\vee}$ be a hyperelliptic K3, i.e.\ a double
cover of the plane ramified over a sextic $\delta$. We will assume
that $\delta$ is generic; in particular it does not admit a tritangent
(this is a codimension one condition on the space of plane
sextics). The pull-back of a generic line in $(\P^2)^{\vee}$ gives a
smooth genus two curve $C$ in $S$, whose linear system is the $\P^2$
dual to $(\P^2)^{\vee}$. We also take $H=C$ as an ample divisor on
$S$. Let us use $Z^2$ to denote the moduli space $Z=\M^s_H((0,[C],1))$
described in the previous subsection. Thus $Z^2$ is a four-fold fibred
by abelian surfaces.
$$\begin{array}{ccc}
\mathrm{Pic}^2 & \hookrightarrow & Z^2 \\
               &                 & \downarrow \\
               &                 & |C|\cong\P^2 \\
\end{array}$$
For all $d\in\Z$, we can construct a similar fibration $Z^d$, whose 
generic fibre is the degree $d$ Picard group $\mathrm{Pic}^dD$ of a
smooth genus two curve $D\in |C|$.
$$\begin{array}{ccc}
\mathrm{Pic}^d & \hookrightarrow & Z^d \\
               &                 & \downarrow \\
               &                 & |C|\cong\P^2 \\
\end{array}$$
Indeed, $Z^d$ is just the Mukai moduli space
$\M^s_H((0,[C],k-1))$. Since the Mukai vectors are primitive, $Z^d$ 
are smooth compact irreducible holomorphic symplectic four-folds,
deformation equivalent to $\mathrm{Hilb}^2(S)$.

As in the first remark in Section~\ref{identify}, $Z^d$ can be
identified with the compactified relative Jacobian
$\overline{\mathrm{Pic}}^d(\mathcal{C}/|C|)$ of the family of curves
$\mathcal{C}\rightarrow |C|$, and 
Markushevich~\cite{markushevich95,markushevich96} gave explicit
constructions of $Z^0$ and $Z^1$ via this approach (we will see
shortly that the other spaces $Z^d$ are all isomorphic to one of these
two). Note that the curves in $|C|$ are
\begin{enumerate}
\item smooth genus two curves, generically; pull-backs of lines in
  $(\P^2)^{\vee}$ meeting $\delta$ transversely,
\item genus one curves with one node, in codimension one; 
pull-backs of lines tangent to $\delta$ at precisely one point,
\item genus one curves with one cusp, in codimension two; 
pull-backs of flex lines of $\delta$,
\item and rational curves with two nodes, in codimension two;
  pull-backs of bitangents to $\delta$.
\end{enumerate}
The singular curves sit above the curve $\Delta\subset\P^2$ dual to
$\delta$. By the Pl{\"u}cker formulae~\cite{gh78} $\Delta$ is a degree
30 curve with 72 cusps and 324 nodes. The two kinds of most singular
fibres sit above the cusps and nodes respectively.

We see in particular that all the curves are integral (i.e.\ reduced
and irreducible). Altman, Iarrobino, and Kleiman~\cite{aik77} proved
that for families of integral curves embedded in surfaces, the fibres
of $\overline{\mathrm{Pic}}^d(\mathcal{C}/|C|)$ are always
irreducible. In fact, this follows from the weaker result of
D'Souza~\cite{dsouza79} since our curves have at worst nodes or cusps
as singularities; D'Souza also shows that the fibres are reduced and
equidimensional in this case. In fact we will give an explicit 
description of all the fibres.

Clearly for a smooth genus two curve $D$, of type (1), the Jacobian 
$\mathrm{Pic}^dD$ is already compact, and is a smooth abelian 
surface. For type (2) we have the following description.

\begin{lemma}
\label{type2}
Let $D$ be a genus one curve with one node $r$, i.e.\ of arithmetic 
genus two. Let $\pi:\tilde{D}\rightarrow D$ be the normalization of 
$D$, and let $p$ and $q$ be the two points of $\tilde{D}$ which are 
identified by $\pi$. Let $\L$ be the Poincar{\'e} line bundle on 
$\tilde{D}\times\mathrm{Pic}^0\tilde{D}$, and let $\L_p$ and $\L_q$ 
be the restrictions to $\{p\}\times\mathrm{Pic}^0\tilde{D}$ and 
$\{q\}\times\mathrm{Pic}^0\tilde{D}$ respectively. Then 
$\overline{\mathrm{Pic}}^0D$ is given by taking the $\P^1$-bundle 
$\P(\L_p\oplus\L_q)$ over $\mathrm{Pic}^0\tilde{D}$ and identifying 
$s_0:=\P(\L_p)\cong\mathrm{Pic}^0\tilde{D}$ and 
$s_{\infty}:=\P(\L_q)\cong\mathrm{Pic}^0\tilde{D}$ with a translation
by $\O(p-q)$. (Note that $\O(p-q)$ is a point on 
$\mathrm{Pic}^0\tilde{D}$, by which we can translate.)
\end{lemma}

\begin{proof}
This is Example (1) on page 83 of~\cite{os79}. First let us make 
one general comment: a torsion-free rank one sheaf on a curve
must be locally free at all smooth points, and at a singular point
$r$ it is either locally free or isomorphic to the maximal ideal
$\mathfrak{m}_r$.

Consider first the locally free case. A line bundle on $D$ is given 
by a line bundle on $\tilde{D}$ plus `gluing data': an element of $\C^*$ 
describing how the fibres at $p$ and $q$ are to be identified. Thus 
we have the extension
$$0\rightarrow{\mathbb G}_m\rightarrow\mathrm{Pic}^0D\rightarrow\mathrm{Pic}^0\tilde{D}\rightarrow 0$$
as in pages 247-253 of~\cite{hm98}.

Another description, which makes it clearer how to compactify, is
to begin with a degree zero line bundle $L$ on $\tilde{D}$ and 
push-forward by $\pi$. The resulting sheaf $\pi_*L$ has fibre 
$L_p\oplus L_q$ at the node $r$. Let $L_0$ be a one-dimensional 
subspace of $L_p\oplus L_q$ and let $\xi$ be the composition
$$\pi_*L\stackrel{\mathrm{ev}_r}{\longrightarrow}L_p\oplus L_q\rightarrow L_p\oplus L_q/L_0$$
where here we regard $L_p$, $L_q$, and $L_0$ as skyscraper 
sheaves supported at $r$. Then $\mathrm{ker}\xi$ is a torsion-free 
rank one degree zero sheaf on $D$.

The line bundle on $\mathrm{Pic}^0\tilde{D}$ whose fibre over $L$ is
$L_p\otimes L_q^{\vee}$ is precisely $\L_p\otimes \L_q^{\vee}$ and
hence the one-dimensional subspace $L_0$ is really a point in
the fibre of $\P(\L_p\oplus \L_q)$ over $L$. The above description
therefore leads to the $\P^1$-bundle $\P(\L_p\oplus \L_q)$ over
$\mathrm{Pic}^0\tilde{D}$. This is not quite the compactified
Jacobian of $D$, but rather the scheme representing the presentation 
functor~\cite{ak90}; in fact, it is the normalization
$\widetilde{\mathrm{Pic}}^0D:=\widetilde{\overline{\mathrm{Pic}}^0D}$ 
of the compactified Jacobian of $D$, as we now explain.

If $L_0=L_p$, then the above composition becomes
$$\pi_*L\rightarrow L_q$$
and $\mathrm{ker}\xi =\pi_*(L(-q))$. For another line bundle 
$L^{\prime}$ on $\tilde{D}$, and $L^{\prime}_0=L^{\prime}_q$, we get
$\pi_*(L^{\prime}(-p))$; this is clearly isomorphic to $\pi_*(L(-q))$ 
when $L^{\prime}\cong L\otimes\O(p-q)$. These sheaves represent the 
zero and infinity sections $s_0$ and $s_{\infty}$, and we see that 
they are glued with a translation by $\O(p-q)$.
\end{proof}

Next consider type (3) curves.

\begin{lemma}
\label{type3}
Let $D$ be a genus one curve with one cusp $r$, i.e.\ of arithmetic 
genus two. Let $\pi:\tilde{D}\rightarrow D$ be the normalization of 
$D$, and let $p\in\tilde{D}$ be the preimage of $r$. Let $\L$ be the 
Poincar{\'e} line bundle on $\tilde{D}\times\mathrm{Pic}^0\tilde{D}$ 
and $\L_p$ the restrictions to $\{p\}\times\mathrm{Pic}^0\tilde{D}$. 
Let $J^1\L_p$ denote the first jet bundle of $\L_p$; there is an exact 
sequence
$$0\rightarrow\Omega^1_{\mathrm{Pic}^0\tilde{D}}\otimes\L_p\rightarrow J^1\L_p\rightarrow\L_p\rightarrow 0.$$
Then $\overline{\mathrm{Pic}}^0D$ is given by taking the $\P^1$-bundle 
$\P(J^1\L_p)$ over $\mathrm{Pic}^0\tilde{D}$ and contracting along
$s_{\infty}:=\P(\Omega^1\otimes\L_p)\cong\mathrm{Pic}^0\tilde{D}$ in
a certain direction, to produce a locus of cusps. Note that the 
contraction is not purely in the direction of the fibres.
\end{lemma}

\begin{proof}
This is Theorem 10 of~\cite{kleiman84} (see also~\cite{ak90}).
In the locally free case, we have an extension~\cite{hm98}
$$0\rightarrow{\mathbb G}_a\rightarrow\mathrm{Pic}^0D\rightarrow\mathrm{Pic}^0\tilde{D}\rightarrow 0$$
which once again splits as a sequence of abelian 
groups~\cite{beauville99}. 

We can also get a description of the compactified Jacobian by 
beginning with a degree zero line bundle $L$ on $\tilde{D}$ and 
pushing forward by $\pi$. The resulting sheaf $\pi_*L$ has fibre 
$J^1L_p$ at the node $r$. Let $L_0$ be a one-dimensional subspace 
of $J^1L_p$ and let $\xi$ be the composition
$$\pi_*L\stackrel{\mathrm{ev}_r}{\longrightarrow}J^1L_p\rightarrow J^1L_p/L_0.$$
Then $\mathrm{ker}\xi$ is a torsion-free rank one degree zero 
sheaf on $D$. This leads to the description of the normalization
$\widetilde{\mathrm{Pic}}^0D$ of the compactified Jacobian of $D$ 
as the $\P^1$-bundle $\P(J^1\L_p)$ over $\mathrm{Pic}^0\tilde{D}$. 
As with the previous lemma, some further identifications need to be
made, and these amount to contracting along the locus $s_{\infty}$
in a certain direction. The precise direction is described in 
Kleiman's paper~\cite{kleiman84}.
\end{proof}

Finally, consider type (4) curves.

\begin{lemma}
\label{type4}
Let $D$ be a rational curve with two nodes $r_1$ and $r_2$, i.e.\ of 
arithmetic genus two. Let $\pi:\tilde{D}\rightarrow D$ be the 
normalization of $D$, and let $\{p_1,q_1\}$ and $\{p_2,q_2\}$ be the 
pairs of points of $\tilde{D}$ which are identified by $\pi$. Since 
$\tilde{D}\cong\P^1$, we can define $\lambda$ to be the cross-ratio
of the four points $\{p_1,q_1,p_2,q_2\}$. Note that multiplication by 
$\lambda$ gives an isomorphism $\P^1\rightarrow\P^1$ which fixes $0$
and $\infty$. Then 
$\overline{\mathrm{Pic}}^0D$ is given by taking $\P^1\times\P^1$ and
identifying $s_0^1:=\{0\}\times\P^1$ and 
$s_{\infty}^1:=\{\infty\}\times\P^1$ via multiplication by $\lambda$,
and
identifying $s_0^2:=\P^1\times\{0\}$ and 
$s_{\infty}^2:=\P^1\times\{\infty\}$ via multiplication by $\lambda$.
\end{lemma}

\begin{proof}
This is Example (2) on pages 83-84 of~\cite{os79}. A degree zero
line bundle on $\tilde{D}\cong\P^1$ is necessarily isomorphic to the
trivial bundle. In the locally free case, the gluing data is a point 
in $\C^*\times\C^*$, which tells us how to identify the fibres at 
$p_1$ and $q_1$, and those at $p_2$ and $q_2$, to get a line bundle 
on $D$. Thus
$$\mathrm{Pic}^0D\cong{\mathbb G}_m\times{\mathbb G}_m.$$

To see how to compactify this, begin with a degree zero line bundle 
$L$ on $\tilde{D}$ (isomorphic to the trivial bundle) and 
push-forward by $\pi$. The resulting sheaf $\pi_*L$ has fibres 
$L_{p_1}\oplus L_{q_1}$ and $L_{p_2}\oplus L_{q_2}$ at the nodes $r_1$
and $r_2$, respectively. Let $L_{01}$ and $L_{02}$ be one-dimensional 
subspaces of $L_{p_1}\oplus L_{q_1}$ and $L_{p_2}\oplus L_{q_2}$, 
respectively, and let $\xi$ be the composition
$$\pi_*L\stackrel{\mathrm{ev}_{r_1}\oplus\mathrm{ev}_{r_2}}{\longrightarrow}(L_{p_1}\oplus L_{q_1})\oplus(L_{p_2}\oplus L_{q_2})\rightarrow (L_{p_1}\oplus L_{q_1}/L_{01})\oplus(L_{p_2}\oplus L_{q_2}/L_{02}).$$
Then $\mathrm{ker}\xi$ is a torsion-free rank one degree zero 
sheaf on $D$. Thus we see that the normalization 
$\widetilde{\mathrm{Pic}}^0D$ of the compactified Jacobian of $D$ is 
$$\P(L_{p_1}\oplus L_{q_1})\otimes\P(L_{p_2}\oplus L_{q_2})\cong\P^1\times\P^1.$$

Now if $L_{01}=L_{p_1}$ we find $\mathrm{ker}\xi=\pi_*(L(-q_1))$,
and similarly if $L^{\prime}_{01}=L^{\prime}_{q_1}$ then 
$\mathrm{ker}\xi=\pi_*(L^{\prime}(-p_1))$. Clearly $L^{\prime}$ 
is always isomorphic to $L\otimes\O(p_1-q_1)$, as both are degree one 
line bundles on $\tilde{D}\cong\P^1$. However, tensoring with 
$\O(p_1-q_1)$ changes the other $\P^1\cong\P(L_{p_2}\oplus L_{q_2})$ 
factor. More specifically, if we let $f$ be a meromorphic function on 
$\P^1$ vanishing at $p_1$ and with a simple pole at $q_1$, so that 
$p_1-q_1=(f)$, then
\begin{eqnarray*}
L^{\prime}_{02} & \subset & L^{\prime}_{p_2}\oplus L^{\prime}_{q_2} \\
                & = & (L\otimes\O((f)))_{p_2}\oplus (L\otimes\O((f)))_{q_2}. \\
\end{eqnarray*}
Hence if $L^{\prime}_{02}$ corresponds to $w\in\P^1$, then $L_{02}$
corresponds to $\lambda w\in\P^1$, where $\lambda=f(p_2)/f(q_2)$ is
the cross-ratio of the four points $\{p_1,q_1,p_2,q_2\}$. The same 
factor arises when we glue $s^2_0$ to $s^2_{\infty}$. This completes 
the proof.
\end{proof}

\begin{remark}
\label{all_degrees}
All of the compactified Jacobians described above are irreducible 
(as required by~\cite{aik77}). Therefore the isomorphism
$$\mathrm{Pic}^0D\rightarrow\mathrm{Pic}^dD$$
given by tensoring with some fixed degree $d$ line bundle extends to
an isomorphism
$$\overline{\mathrm{Pic}}^0D\rightarrow\overline{\mathrm{Pic}}^dD$$
and we have a description of the compactified Picard schemes of all
degrees.
\end{remark}

\begin{remark}
Only the compactified Jacobian of type (4) has non-zero Euler 
characteristic, equal to one. Therefore only these fibres make a
non-trivial contribution to the Euler characteristic of $Z^d$, which
is therefore 324 (see Beauville~\cite{beauville99} for how this
method can be used to calculate the Euler characteristic of
$\mathrm{Hilb}^nS$).
\end{remark}

\begin{remark}
\label{flat}
A fact that we will use later is that $Z^d\rightarrow\P^2$ is a flat
fibration. This follows from the corollary after Theorem 23.1 in
Matsumura's book~\cite{matsumura86}, since $Z^d$ is smooth and the
fibres are equidimensional.
\end{remark}





We can use the first of the above remarks to prove the next lemma.


\begin{lemma}
\label{after_flat}
All the spaces $Z^{2m}$ (for $m\in\Z$) are isomorphic and all the 
spaces $Z^{2m+1}$ (for $m\in\Z$) are isomorphic, so we essentially 
just have $Z^0$ and $Z^1$. Moreover, $Z^1$ is a torsor over $Z^0$.
\end{lemma}

\begin{proof}
First note that $Z^0$ admits a global section, given by taking the
trivial (degree zero) line bundle $\O_D$ on each curve $D\in
|C|$. Since the fibres of $Z^d$ are reduced and irreducible, choosing
a local section of $\mathcal{C}$ over $U\subset |C|$ will give a local
isomorphism $Z^d/U\rightarrow Z^{d+1}/U$, and hence all the spaces
$Z^d$ are locally isomorphic as fibrations. In particular, this means
$Z^d$ is a torsor over $Z^0$.

Now for all $d\in\Z$ there is a global isomorphism 
$$Z^d\stackrel{\cong}{\longrightarrow}Z^{d+2}$$
which over $D\in|C|$ is given by tensoring stable sheaves with the
canonical bundle $\K_D$ (of degree two). This completes the proof.
\end{proof}

\begin{proposition}
For a generic hyperelliptic K3 surface $S$ (i.e.\ the sextic $\delta$
is generic), the spaces $Z^0$ and $Z^1$ are not isomorphic. Indeed
they have different periods, and hence are not even birational.
\end{proposition}

\begin{proof}
We will use O'Grady's description~\cite{ogrady97} of the weight two 
Hodge structure of the Mukai moduli space $\M^s_H(v)$ to show that 
$Z^0$ and $Z^1$ have non-isomorphic Picard lattices, and hence 
different periods. The Mukai lattice $\H^{\bullet}(S,\Z)$ can be 
given the Hodge structure whose $(2,0)$, $(1,1)$, and $(0,2)$ 
components are
$$\H^{2,0}(S),\qquad\qquad\H^{0,0}(S)\oplus\H^{1,1}(S)\oplus\H^{2,2}(S),\qquad\mbox{and}\qquad\H^{0,2}(S)$$
respectively. Then O'Grady proved that for $(v,v)>0$, the weight two
Hodge structure of $\M^s_H(v)$ is isomorphic to $v^{\perp}$ (and the
induced quadratic form agrees with the Beauville-Bogomolov quadratic
form on $\H^2(\M^s_H(v),\Z)$). For generic $S$ we can assume that the
Picard lattice $\H^{1,1}(S)\cap\H^2(S,\Z)$ is generated by $[C]$. Then
the Picard lattice of $\M^s_H(v)$ is isomorphic to
$$\{(a,b,c)\in\Z^3|(a,b[C],c)\in v^{\perp}\}.$$

For the particular spaces we are interested in, we find the Picard 
lattices of $Z^0=\M^s_H((0,[C],-1))$ and $Z^1=\M^s_H((0,[C],0))$ are
$$\Z(-2,[C],0)\oplus\Z(0,0,1)$$
and
$$\Z(-1,0,0)\oplus\Z(0,0,1)$$
respectively, and the induced quadratic forms
$$\left(\begin{array}{cc}
        2 & 2 \\
        2 & 0 \\
        \end{array}\right)\qquad\mbox{and}\qquad
\left(\begin{array}{cc}
        0 & 1 \\
        1 & 0 \\
        \end{array}\right)$$
are easily seen to be non-isomorphic (only the second is unimodular).
\end{proof}

\begin{corollary}
Generically, $Z^1$ does not admit a section. It does, however, admit a
rational 2-valued section.
\end{corollary}

\begin{proof}
If $Z^1$ admitted a section it would be isomorphic to $Z^0$. A
2-valued section is a pair of degree one line bundles on each curve
$D\in |C|$ (since we will construct a rational 2-valued section, we
can avoid the singular fibres, and therefore talk about genuine line
bundles rather than stable rank one sheaves). Note that a degree two
line bundle (such as $\K_D$) on each curve does {\em not\/} 
canonically give a pair of degree one line bundles. However, by 
adjunction
$$\K_D=\O(D)|_D=\O(C)|_D$$
as $\K_S$ is trivial. In other words, $C|_D$ is just a pair of points
on each curve $D$, giving us a pair of degree one line bundles as
required. 
\end{proof}

\begin{remark}
This is only a rational 2-valued section as obviously $C$ does not
intersect itself transversely. A genuine 2-valued section appears not
to exist, though this will not present any difficulties. From the
discussion in Section 2, we know that there should exist an element
$\alpha\in\H^1(\P^2,Z^0)$ classifying the torsor $Z^1$. The fact that
$Z^1$ admits a rational 2-valued section means that $\alpha$ should be
2-torsion, at least in {\'e}tale cohomology
$\H^1_{\mathrm{\acute{e}t}}(\P^2,Z^0)$. This element was constructed 
explicitly by Markushevich in~\cite{markushevich96}.
\end{remark}

\begin{remark}
An example of a non-generic hyperelliptic K3 surface $S$ arises when
the sextic admits a tritangent: this condition is codimension one on
the space of sextics. In this case the pull-back of the tritangent to
$S$ gives a reducible curve, consisting of two rational curves $C_1$
and $C_2$ which intersect transversely at three points. Let us assume
that $S$ is otherwise as generic as possible, so that its Picard
lattice is generated by $[C_1]$ and $[C_2]$ (since $[C]=[C_1]+[C_2]$
and $C.C=2$, $[C_1]$ and $[C_2]$ cannot be proportional). The Picard
lattice of $\M^s_H(v)$ is therefore isomorphic to
$$\{(a,b,c,d)\in\Z^4|(a,b[C_1]+c[C_2],d)\in v^{\perp}\}.$$

A bit of work shows that now $Z^0$ and $Z^1$ have isomorphic Picard
lattices, with quadratic form
$$\left(\begin{array}{ccc}
        0 & 1 & 0 \\
        1 & 0 & 0 \\
    0 & 0 & -10 \\
        \end{array}\right)$$
with respect to the bases
$$\Z(-1,[C_1],1)\oplus\Z(-5,4[C_1]+[C_2],1)\oplus\Z(-10,9[C_1]+[C_2],5)$$
and
$$\Z(-1,0,0)\oplus\Z(0,0,1)\oplus\Z(0,[C_1]-[C_2],0)$$
respectively.

In this case, there is a birational map $Z^0\dashrightarrow Z^1$
inducing the isomorphism of periods: since $C_1.C=1$ the restriction
of the line bundle $\O(C_1)$ on $S$ to a generic curve $D\in |C|$ is a
degree one line bundle, and tensoring with this line bundle induces
the birational map. This line bundle also provides a (birational)
section of $Z^1$.
\end{remark}

\section{Derived equivalences}

In Section 3 we described a certain four-fold $Z^0$ which is fibred 
by abelian surfaces, as well as another fibration $Z^1$ which is a 
torsor over $Z^0$. In this section we will reinterpret this torsor 
via moduli spaces and twisted Fourier-Mukai transforms. In 
particular, we will show that the derived category of sheaves on 
$Z^1$ is equivalent to the derived category of twisted sheaves on 
$Z^0$.

\subsection{The dual fibration}

Since $Z^1$ is fibred over $\P^2$, we can define its compactified 
relative Picard scheme, which we regard as a dual fibration. Our
first goal is to understand this space. Firstly, observe that $Z^0$ 
and $Z^1$ are locally isomorphic as fibrations, and therefore their
dual fibrations $\overline{\mathrm{Pic}}^0(Z^0/\P^2)$ and 
$\overline{\mathrm{Pic}}^0(Z^1/\P^2)$ must also be locally 
isomorphic. Since these dual fibrations both admit canonical global 
sections (given by the trivial bundle on each fibre of $Z^0$, 
respectively $Z^1$), we must have
$$\overline{\mathrm{Pic}}^0(Z^0/\P^2)\cong\overline{\mathrm{Pic}}^0(Z^1/\P^2).$$
In this section we will define $P$ to be 
$\overline{\mathrm{Pic}}^0(Z^0/\P^2)$, and we will show that $P$ is 
also isomorphic to $Z^0$; thus $Z^0$ is self-dual. By the above
observations, this also implies that $Z^0$ is the dual fibration
of $Z^1$ (the operation of taking the dual fibration is only locally 
reflexive). Thus we will interpret $Z^0$, which is a priori a moduli 
space of sheaves on the K3 surface $S$, as also being a moduli space 
of stable sheaves on the four-fold $Z^1$.

It is clear that corresponding smooth fibres of $Z^0$ and $P$ are
isomorphic: if $\mathrm{Pic}^0D$ is a smooth fibre of $Z^0$ (i.e.\ 
$D$ is of type (1), a smooth genus two curve), then the corresponding 
fibre of $P$ is
$$\mathrm{Pic}^0(\mathrm{Pic}^0D)\cong\mathrm{Pic}^0D.$$
This result also extends to the singular fibres by the 
following autoduality result of Esteves and Kleiman. Firstly, 
since the fibres of $Z^0$ are reduced and irreducible, recall that 
by Lemma~\ref{no_fat_sheaves} the fibres of $P$ are just the 
compactified Picard schemes of the fibres of $Z^0$. 

\begin{theorem}[Esteves and Kleiman~\cite{ek04}]
\label{autoduality}
Let $D$ be a surficial curve with at worst nodes and/or cusps as
singularities. Then the compactified Picard scheme of the 
compactified Jacobian of $D$ is isomorphic to the compactified 
Jacobian of $D$, i.e.\
$$\overline{\mathrm{Pic}}^0(\overline{\mathrm{Pic}}^0D)\cong\overline{\mathrm{Pic}}^0D.$$
\end{theorem}

\begin{remark}
In an earlier paper Esteves, Gagn{\'e}, and Kleiman~\cite{egk02}
proved that if $D$ is a surficial curve with at worst double points as
singularities (which includes nodes or cusps), then
$$\mathrm{Pic}^0(\overline{\mathrm{Pic}}^0D)\cong\mathrm{Pic}^0D.$$
Now it is not immediately true that the same result must hold for 
compactified Picard schemes. A priori, a compactified Picard scheme
could consist of several irreducible components, with whole 
components parametrizing torsion-free rank one sheaves {\em none\/} 
of which are locally free. Such examples even exist (see Altman, 
Iarrobino, and Kleiman~\cite{aik77} for an example of a space curve 
whose compactified Jacobian has this property), but only in the 
presence of particularly bad kinds of singularities. It was also
shown in~\cite{aik77} that the compactified Jacobian of a surficial
curve must be irreducible. So when the singularities are mild, as in
the theorem, one expects much better behaviour.
\end{remark}

In particular, Theorem~\ref{autoduality} applies to our curves of type 
(2), (3), and (4). Let us take a closer look at the type (2) case.

\begin{lemma}
If $D$ is a curve of type (2), then the compactified Picard scheme of 
the compactified Jacobian of $D$ is isomorphic to the compactified 
Jacobian of $D$, i.e.\
$$\overline{\mathrm{Pic}}^0(\overline{\mathrm{Pic}}^0D)\cong\overline{\mathrm{Pic}}^0D.$$
\end{lemma}

\begin{proof}
Our argument will identify the two spaces {\em as sets}. More care
needs to be taken to prove that their analytic structures agree (for
instance, one could exhibit a universal sheaf for the moduli problem);
for this we refer to Esteves and Kleiman's proof of the more general
Theorem~\ref{autoduality}.

We begin by simplifying the description
of $\overline{\mathrm{Pic}}^0D$ given in Lemma~\ref{type2}. Pick an
arbitrary point $p_0$ in the normalization $\tilde{D}$. Then we have
an isomorphism $\tilde{D}\cong\mathrm{Pic}^0\tilde{D}$ given by 
taking the point $p\in\tilde{D}$ to $\O(p-p_0)$. Under this 
isomorphism the bundle $\L_p=\L|_{\{p\}\times\mathrm{Pic}^0\tilde{D}}$ 
on $\mathrm{Pic}^0\tilde{D}$ becomes the bundle
$$\L|_{\tilde{D}\times\{\O(p-p_0)\}}=\O(p-p_0)$$
on $\tilde{D}$. Translation by 
$$\O(p-q)=\O(p-p_0)\otimes\O(q-p_0)^{\vee}$$
on $\mathrm{Pic}^0\tilde{D}$ becomes translation by $p-q$ on 
$\tilde{D}$. Thus the compactified Jacobian of $D$ is isomorphic to
the $\P^1$-bundle $\P(\O(p)\oplus\O(q))$ over $\tilde{D}$ with 
$s_0:=\P(\O(p))\cong\tilde{D}$ and 
$s_{\infty}:=\P(\O(q))\cong\tilde{D}$ identified with a
translation by $p-q$. The normalization $\widetilde{\mathrm{Pic}}^0D$ 
of the compactified Jacobian is isomorphic to $\P(\O(p)\oplus\O(q))$ 
itself.

Now consider $\mathrm{Pic}^0(\overline{\mathrm{Pic}}^0D)$. We will 
begin with a line bundle $\E$ on the normalization 
$\widetilde{\mathrm{Pic}}^0D$, and then glue $\E|_{s_0}$ to 
$\E|_{s_{\infty}}$ to obtain a sheaf on 
$\overline{\mathrm{Pic}}^0D$. We choose $\E$ to be the pull-back 
$\gamma^*L$ of a degree zero line bundle $L$ on $\tilde{D}$, under 
the projection
$$\gamma:\P(\O(p)\oplus\O(q))\rightarrow\tilde{D}.$$
Both $s_0$ and $s_{\infty}$ are identified with $\tilde{D}$, by the 
map $\gamma$; therefore
$$\E|_{s_0}=(\gamma|_{s_0})^*L$$
and
$$\E|_{s_{\infty}}=(\gamma|_{s_{\infty}})^*L.$$
However, $s_0$ and $s_{\infty}$ are not glued by 
$(\gamma|_{s_{\infty}})^{-1}\circ(\gamma|_{s_0})$, but rather by
$$\tau:=(\gamma|_{s_{\infty}})^{-1}\circ\mathrm{tr}\circ(\gamma|_{s_0}):s_0\rightarrow s_{\infty}$$
where $\mathrm{tr}$ is translation by $p-q$. Any degree zero line 
bundle on $\tilde{D}$ can be written as $\O(a-b)$ for some points 
$a$ and $b\in\tilde{D}$; therefore 
$$\mathrm{tr}^*\O(a-b)\cong\O((a+p-q)-(b+p-q))\cong\O(a-b).$$
In particular, $\mathrm{tr}^*L\cong L$, and so we still have
$$\tau^*(\E|_{s_{\infty}})\cong\E|_{s_0}.$$
Choosing a non-zero element of
$$\mathrm{Hom}_{s_0}(\E|_{s_0},\tau^*(\E|_{s_{\infty}}))\cong\mathrm{Hom}_{\tilde{D}}(L,\mathrm{tr}^*L)\cong\C$$
thus enables us to glue $\E|_{s_0}$ to $\E|_{s_{\infty}}$, giving a
rank-one locally free sheaf on $\overline{\mathrm{Pic}}^0D$. Thus we
have described the Picard scheme of $\overline{\mathrm{Pic}}^0D$ as a 
$\C^*$-bundle over $\mathrm{Pic}^0\tilde{D}$. It will become clear 
when we compactify that this is, in fact, the same as the 
$\C^*$-bundle arising in the description of $\mathrm{Pic}^0D$, and 
thus
$$\mathrm{Pic}^0(\overline{\mathrm{Pic}}^0D)\cong\mathrm{Pic}^0D$$
as required by~\cite{egk02}.

To compactify, we imitate the proofs of Lemmas~\ref{type2}, 
\ref{type3}, and \ref{type4}. First push $\E$ forward under the 
normalization map
$$\pi:\widetilde{\mathrm{Pic}}^0D\rightarrow\overline{\mathrm{Pic}}^0D.$$
Over the singular locus $s\cong\tilde{D}$, $\E$ looks like
$$\E|_{s_0}\oplus\tau^*(\E|_{s_{\infty}})$$
on $s_0$, or equivalently, $L\oplus\mathrm{tr}^*L$ on $\tilde{D}$.
Choose a degree zero line subbundle $L_0$ of $L\oplus\mathrm{tr}^*L$.
Let $\xi$ be the composition
$$\pi_*\E\stackrel{|_{s}}{\rightarrow}L\oplus\mathrm{tr}^*L\rightarrow 
L\oplus\mathrm{tr}^*L/L_0$$
where we regard $L\oplus\mathrm{tr}^*L$ and $L_0$ as torsion sheaves 
supported on $s$. Then $\mathrm{ker}\xi$ is a torsion-free rank one 
sheaf on $\overline{\mathrm{Pic}}^0D$, which is clearly in the same
connected component as the trivial sheaf.

Since $\mathrm{tr}^*L\cong L$, we must have $L_0\cong L$. Thus the 
line subbundle $L_0$ is completely determined by examining the fibres 
over a single point in $\tilde{D}$. Choosing this point to be $p$, we
see that
$$(L_0)_p\subset (L\oplus\mathrm{tr}^*L)_p=L_p\oplus L_q$$
is given by a point of $\P(L_p\oplus L_q)$. Therefore we arrive at a
description of the normalization 
$$\widetilde{\mathrm{Pic}}^0(\overline{\mathrm{Pic}}^0D):=\widetilde{\overline{\mathrm{Pic}}^0(\overline{\mathrm{Pic}}^0D)}$$
of the compactified Jacobian of $\overline{\mathrm{Pic}}^0D$, namely
it is the total space of the $\P^1$-bundle $\P(\L_p\oplus\L_q)$ over
$\mathrm{Pic}^0\tilde{D}$. In particular, comparing to 
Lemma~\ref{type2} we see that
$$\widetilde{\mathrm{Pic}}^0(\overline{\mathrm{Pic}}^0D)\cong\widetilde{\mathrm{Pic}}^0D.$$

Finally, we need to identify some points in 
$\widetilde{\mathrm{Pic}}^0(\overline{\mathrm{Pic}}^0D)$. If we
choose $(L_0)_p$ to be $L_p$ then we find
$$\mathrm{ker}\xi=\pi_*(\E(-s_{\infty}))=\pi_*((\gamma^*L)(-s_{\infty})).$$
For another sheaf $\E^{\prime}=\gamma^*L^{\prime}$, with the choice
$(L^{\prime}_0)_p=L^{\prime}_q$ we get
$$\pi_*(\E^{\prime}(-s_0))=\pi_*((\gamma^*L^{\prime})(-s_0)).$$
Let $f_p$ and $f_q$ be sections of $\O(p)$ and $\O(q)$, respectively, 
on $\tilde{D}$. Define a map
$$\P(\O(p)\oplus\O(q))\rightarrow\P(\O(p+q)\oplus\O(q+p))\cong\P(\O\oplus \O)$$
by taking the point $[x,y]$ over $t\in\tilde{D}$ to 
$[xf_q(t),yf_p(t)]$. Since $\P(\O\oplus\O)$ is the trivial 
$\P^1$-bundle on $\tilde{D}$, this defines a 
meromorphic function on $\P(\O(p)\oplus\O(q))$. Moreover, it has a 
zero when $x$ or $f_q$ vanish (i.e.\ along $s_0$ and $\gamma^{-1}(q)$)
and a pole when $y$ of $f_p$ vanish (i.e.\ along $s_{\infty}$ and
$\gamma^{-1}(p)$). Therefore we have a linear equivalence of 
divisors
$$s_0+\gamma^{-1}(q)\sim s_{\infty}+\gamma^{-1}(p).$$
It follows that the sheaves $(\gamma^*L^{\prime})(-s_0)$ and
$(\gamma^*L)(-s_{\infty})$ are isomorphic when 
$L^{\prime}\cong L\otimes\O(p-q)$. This is precisely the same 
translation that occurred for the compactified Jacobian of $D$.

We have made numerous choices, so it is not clear a priori that we
have a complete description of the compactified Picard scheme of
$\overline{\mathrm{Pic}}^0D$, which by definition is the connected
component of the trivial sheaf in the moduli space of torsion-free
rank one sheaves on $\overline{\mathrm{Pic}}^0D$. At a point which 
represents a locally free sheaf, one can compute the dimension of 
this moduli space and show that it is smooth. Therefore another
irreducible component could only meet the one we have just 
described along its singular locus.

To construct other irreducible components, we could have 
started with a different torsion-free rank one sheaf $\E$ on
$\widetilde{\mathrm{Pic}}^0D$, which in general would look like a
line bundle tensored with the ideal sheaf of a zero-dimensional
subscheme $Z\subset\widetilde{\mathrm{Pic}}^0D$. Since $s_0$ and 
$s_{\infty}$ are smooth curves, the restrictions $\E|_{s_0}$ and 
$\E|_{s_{\infty}}$ would once again be line bundles, and so the 
rest of our construction would still make sense. We could also
have chosen the line subbundle $L_0$ differently; for instance, of 
degree less than zero. All of these variations produce other 
irreducible components of the moduli space of torsion-free rank 
one sheaves, but one can check directly that they don't meet the
irreducible component that we first described. For instance, the 
Euler characteristic of a sheaf must remain constant in a
connected component, and therefore it can be used to distinguish
many connected components (though this on its own is by no means
sufficient). Thus our description
$$\overline{\mathrm{Pic}}^0(\overline{\mathrm{Pic}}^0D)\cong\overline{\mathrm{Pic}}^0D$$
for type (2) curves is complete.
\end{proof}

\begin{remark}
The above argument can probably be extended to type (3) and (4) 
curves: we leave this as an exercise for the (patient) reader.
\end{remark}

\begin{corollary}
\label{PisZ0}
The spaces $P$ and $Z^0$ are isomorphic. In particular, $P$ is 
a smooth holomorphic symplectic four-fold.
\end{corollary}

\begin{proof}
Because of Theorem~\ref{autoduality}, we know that
$$\overline{\mathrm{Pic}}^0(\overline{\mathrm{Pic}}^0D)\cong\overline{\mathrm{Pic}}^0D$$
for all four kind of curves. Moreover, Esteves and Kleiman's result
(Theorem 4.1 in~\cite{ek04}) also holds in the relative case. So for
the family of curves $\mathcal{C}\rightarrow\P^2$, we have
$$\overline{\mathrm{Pic}}^0(\overline{\mathrm{Pic}}^0(\mathcal{C}/\P^2))\cong\overline{\mathrm{Pic}}^0(\mathcal{C}/\P^2)$$
which is precisely $P\cong Z^0$.
\end{proof}

\begin{remark}
\label{PisZ0remark}
Even if we didn't have the relative version of 
Theorem~\ref{autoduality}, we could still conclude that $Z^0\cong P$ 
by arguing that they are locally isomorphic as fibrations and then
observing that they both admit global sections. To prove they are
locally isomorphic we can take a (Stein) open cover $\{U_i\}$ of 
$\P^2$ and consider the induced maps of each open set $U_i$ into a 
compactification of the moduli space of principally polarized 
abelian surfaces. Then over $U_i$, both $Z^0$ and $P$ will be 
given by pulling back the same universal variety.

The difficulty with this argument is that there are different
compactifications of the moduli space of abelian surfaces, not all
of which admit a universal variety. If we use Mumford's
compactification, then there exist local universal varieties 
(see Hulek, Kahn, and Weintraub~\cite{hkw93}), which suffice since 
the $U_i$ are Stein open sets. However, not all degenerations of 
abelian surfaces arise in Mumford's compactification: we get 
compactified Jacobians of type (2) curves, but not of types (3) 
and (4). To proceed we first remove from $U_i$ the points 
corresponding to type (3) and (4) curves 
(fortunately $U_i$ remains simply-connected). We conclude then that
$Z^0$ and $P$ are isomorphic as fibrations over the complement of
a finite set of points in $\P^2$. This also gives a birational map
between $Z^0$ and $P$. Since this map is an isomorphism in 
codimension one, the holomorphic symplectic form $\sigma$ on $Z^0$
induces a holomorphic two-form on $P$ (a priori, we don't know
whether it is non-degenerate).

Now the argument becomes a little delicate: we can show using 
Fourier-Mukai methods (see below) that $P$ is smooth and derived
equivalent to $Z^0$. Then the fact that $Z^0$ is holomorphic
symplectic implies the same is true of $P$ (the derived equivalence 
preserves the Serre functor, so $P$ certainly has trivial canonical 
bundle; this implies the holomorphic two-form on $P$ is indeed
non-degenerate). Finally, a birational map between holomorphic 
symplectic four-folds must be a composition of Mukai flops of 
embedded $\P^2$s (see Burns, Hu, and Luo~\cite{bhl03}), but the 
indeterminacy of our birational map is contained in the (disjoint)
union of the set of singular fibres corresponding to curves of 
types (3) and (4). In particular, it contains no $\P^2$s and thus
the birational map must extend to an isomorphism.
\end{remark}

\subsection{A Fourier-Mukai transform}

We have seen that $P$ parametrizes stable sheaves on $Z^0$. Moreover,
since $Z^0$ admits a global section, there exists a universal sheaf 
$\U$ on $Z^0\times P$. We can therefore construct the functor
$$\Phi^{\U}_{P\rightarrow Z^0}:\D (P)\rightarrow\D (Z^0).$$ 

\begin{theorem}
\label{untwisted_equivalence}
The functor $\Phi^{\U}_{P\rightarrow Z^0}$ is an equivalence of
triangulated categories.
\end{theorem}

\begin{proof}
We will apply Bridgeland and Maciocia's 
Theorem~\ref{bridgeland_maciocia}. Let
$\U_m:=\Phi^{\U}_{P\rightarrow Z^0}\O_m$ be the sheaf on $Z^0$ which
the point $m\in P$ parametrizes. We must show
\begin{enumerate}
\item for all $m\in P$, $\U_m\otimes\K_{Z^1}=\U_m$ and $\U_m$ is
  simple,
\item for all $m_1\neq m_2\in P$,
$$\mathrm{Hom}_{Z^0}(\U_{m_1},\U_{m_2})=0,$$
and the closed subscheme
$$\Gamma(\U):=\{(m_1,m_2)\in P\times
P|\mathrm{Ext}^i_{Z^0}(\U_{m_1},\U_{m_2})\neq 0\mbox{ for some }i\in\Z\}$$
of $P\times P$ has dimension at most five.
\end{enumerate}
Since $\U_m$ is stable, it is also simple. Since $Z^0$ is holomorphic
symplectic, it has trivial canonical bundle $\K_{Z^0}$. Thus condition
(1) follows.

If $\U_{m_1}\rightarrow\U_{m_2}$ is a non-trivial morphism, then it
must be an isomorphism since $\U_{m_2}$ is stable. Therefore
$m_1=m_2$, proving the first part of condition (2). It remains to prove
the bound on the dimension of $\Gamma(\U)$.

Firstly, suppose that $m_1$ and $m_2$ lie in different fibres of
$P$. Then $\U_{m_1}$ and $\U_{m_2}$ are sheaves supported on different
(disjoint) fibres of $Z^0$. Therefore
$$\mathrm{Ext}^i_{Z^0}(\U_{m_1},\U_{m_2})=0$$
for all $i$ because all local $\mathcal{E}xt$ sheaves vanish.

Next suppose that $m_1$ and $m_2$ lie in the same {\em smooth\/} fibre 
of $P$. Then $\U_{m_1}$ and $\U_{m_2}$ are supported on the smooth 
fibre $Z^0_t=\mathrm{Pic}^0D$ of $Z^0$; in fact they are of the form
$\iota_*L_1$ and $\iota_*L_2$ respectively, where
$\iota:Z^0_t\hookrightarrow Z^0$ is inclusion and $L_1$ and $L_2$ are
degree zero line bundles on the (smooth) abelian surface
$Z^0_t$. Now we have the following spectral sequence, taken from
Section 7.2 of Bridgeland and Maciocia~\cite{bm02}
$$E^{p,q}_2:=\mathrm{Ext}^p_{Z^0_t}(L_1\otimes\Lambda^q\O_{Z^0_t}^{\oplus
  2},L_2)\Longrightarrow\mathrm{Ext}^{p+q}_{Z^0}(\iota_*L_1,\iota_*L_2)$$
where $\O_{Z^0_t}^{\oplus 2}$ is really the conormal bundle of $Z^0_t$
in $Z^0$, which is trivial. If $m_1\neq m_2$ then $L_1$ and $L_2$ are
not isomorphic. The cohomology of the non-trivial degree zero line
bundle $L_1^{\vee}\otimes L_2$ on $Z^0_t$ therefore vanishes in all
degrees (see Chapter 3 of Birkenhake and Lange~\cite{bl92}), i.e.\ 
$$\H^p(Z^0_t,L_1^{\vee}\otimes L_2)=0\mbox{ for all }p\in\Z.$$
It follows that the spectral sequence vanishes and we have proved
$$\mathrm{Ext}^i_{Z^0}(\U_1,\U_2)=\mathrm{Ext}^i_{Z^0}(\iota_*L_1,\iota_*L_2)=0$$ 
for all $i\in\Z$ in this case.

We have shown that $\Gamma(\U)$ is a subset of
$$\mathrm{Diag}\cup\{(m_1,m_2)\in P\times P|m_1\mbox{ and }m_2\mbox{ lie 
in the same singular fibre}\}$$
where $\mathrm{Diag}$ is the diagonal in $P\times P$. But the singular
fibres of $P$ sit above the curve $\Delta\subset\P^2$ and have
dimension two. Thus
$$\{(m_1,m_2)\in P\times P|m_1\mbox{ and }m_2\mbox{ lie in
  the same singular fibre}\}$$
has dimension five, and $\mathrm{Diag}$ clearly has dimension
four. Therefore $\Gamma(\U)$ has dimension at most five and the
theorem is proved.
\end{proof}

\begin{remark}
Theorem~\ref{bridgeland_maciocia} also implies that $P$ is smooth. We 
used this fact in the remark following Corollary~\ref{PisZ0}, to give
an alternate proof of the isomorphism $P\cong Z^0$.
\end{remark}


\begin{remark}
Since $P\cong Z^0$ by Corollary~\ref{PisZ0}, the equivalence 
$$\Phi^{\U}_{P\rightarrow Z^0}:\D(P)\rightarrow\D(Z^0)$$
is really an auto-equivalence of the derived category of $Z^0$.
It is non-trivial: for instance, on a smooth fibre it induces the
non-trivial auto-equivalence of the derived category of a 
principally polarized abelian surface (which was first constructed 
by Mukai~\cite{mukai81}).
\end{remark}

In fact, we can extend Mukai's results to degenerate abelian 
surfaces as follows. Let $D$ be an arbitrary curve of type (2), (3), or (4).
We claim that $D$ occurs in some hyperelliptic K3 surface $S$. This
is straightforward: we write $D$ as a double cover of $\P^1$ 
branched over six points, some of them coinciding. Then we embed 
the line $\P^1$ in $(\P^2)^{\vee}$ and find a sextic $\delta$
meeting the line in the given six points, with $\delta$ touching 
the line to order two or three at a point of multiplicity two, 
respectively three. The space of plane sextics moduli automorphisms
of $\P^2$ has dimension 19, whereas the points impose just six 
conditions, so we clearly can find such a sextic, and moreover we
can assume that it does not admit a tritangent, so that the 
corresponding K3 surface $S$ is generic.

Next we construct $Z^0$ and $P$ as before. The compactified Jacobian
$J:=\overline{\mathrm{Jac}}D$ occurs as a fibre of $Z^0$, and its dual  
$\overline{\mathrm{Pic}}^0J$ is the corresponding fibre of 
$P$. Use $\iota$ to denote the inclusion of $J$ in $Z^0$,
$\overline{\mathrm{Pic}}^0J$ in $P$, and
$J\times\overline{\mathrm{Pic}}^0J$ in $Z^0\times P$ (the usage will
always be clear from the context). The universal sheaf $\U$ on
$Z^0\times P$ restricts to a universal sheaf $\U_D:=\iota^*\U$ on
$J\times\overline{\mathrm{Pic}}^0J$, which generalizes the
Poincar{\'e} line bundle for a smooth abelian surface.

\begin{corollary}
\label{singular_Mukai}
Let $D$ be an arbitrary curve of type (2), (3), or (4), and let
$J:=\overline{Jac}D$ be its compactified Jacobian, which is a
degeneration of a principally polarized abelian surface. The 
universal sheaf $\U_D$ induces an integral transform
$$\Phi^{\U_D}_{\overline{\mathrm{Pic}}^0J\rightarrow J}:\D(\overline{\mathrm{Pic}}^0J)\rightarrow \D(J)$$
which is an equivalence of triangulated categories.
\end{corollary}

\begin{proof}
Our proof is taken directly from Section 6 of
Chen~\cite{chen02}. Write $\Phi$ for $\Phi^{\U}_{P\rightarrow Z^0}$
and $\Phi_D$ for $\Phi^{\U_D}_{\overline{\mathrm{Pic}}^0J\rightarrow
J}$. Label the various projections as in the following diagram:
$$\begin{array}{ccccc}
Z^0 & \stackrel{\pi_0}{\leftarrow} & Z^0\times P &
\stackrel{\pi_P}{\rightarrow} & P \\
\uparrow\iota & & \uparrow\iota & & \uparrow\iota \\
J & \stackrel{\pi^D_0}{\leftarrow} & J\times\overline{\mathrm{Pic}}^0J &
\stackrel{\pi^D_P}{\rightarrow} & \overline{\mathrm{Pic}}^0J \\
\end{array}$$

We first claim there is a natural isomorphism of functors
$$\Phi\circ\iota_*\cong\iota_*\circ\Phi_D$$
(this is Lemma 6.1 in~\cite{chen02}). Let
$\E\in\D(\overline{\mathrm{Pic}}^0J)$; using flat base change and the
projection formula we have
\begin{eqnarray*}
\Phi(\iota_*\E) & = & {\bf R}\pi_{0*}(\U\stackrel{\bf L}{\otimes}{\bf L}\pi_P^*(\iota_*\E)) \\
 & \cong & {\bf R}\pi_{0*}(\U\stackrel{\bf L}{\otimes}{\bf
 R}\iota_*({\bf L}(\pi_P^D)^*\E)) \\
 & \cong & {\bf R}\pi_{0*}({\bf
 R}\iota_*(\iota^*\U\stackrel{\bf L}{\otimes}{\bf L}(\pi_P^D)^*\E) \\
 & \cong & {\bf R}\iota_*({\bf
 R}(\pi_0^D)_*(\U_D\stackrel{\bf L}{\otimes}{\bf L}(\pi_P^D)^*\E) \\
 & = & {\bf R}\iota_*(\Phi_D\E) \\
\end{eqnarray*}
and the claim follows.

Now let us prove that $\Phi_D$ is an equivalence (this is Lemma 6.2
in~\cite{chen02}). Let $\Psi$ be the right adjoint of $\Phi$; then
$\Psi$ is an integral transform given by some sheaf $\V$ on $P\times
Z^0$. Let 
$$\Psi_D:\D(J)\rightarrow\D(\overline{\mathrm{Pic}}^0J)$$
be the integral transform given by $\iota^*\V$. Then one can show there
is a natural isomorphism of functors
$$\Psi\circ\iota_*\cong\iota_*\Psi_D$$
using the same argument as for the claim above. We therefore have a
commutative diagram
$$\begin{array}{ccccc}
{\D}(P) & \stackrel{\Phi}{\rightarrow} & {\D}(Z^0) &
\stackrel{\Psi}{\rightarrow} & {\D}(P) \\
\uparrow\iota_* & & \uparrow\iota_* & & \uparrow\iota_* \\
{\D}(\overline{\mathrm{Pic}}^0J) & \stackrel{\Phi_D}{\rightarrow} & {\D}(J) &
\stackrel{\Psi_D}{\rightarrow} & {\D}(\overline{\mathrm{Pic}}^0J). \\
\end{array}$$
In particular
$$(\Psi\circ\Phi)\circ\iota_*\cong\iota_*\circ(\Psi_D\circ\Phi_D),$$
but $\Psi\circ\Phi$ is equivalent to the identity functor on
$\D(P)$. As $\iota$ is a closed embedding, this implies
$\Psi_D\circ\Phi_D$ is equivalent to the identity functor on
$\D(\overline{\mathrm{Pic}}^0J)$. Similarly, using
$\Phi\circ\Psi\cong\mathrm{Id}_{\D(Z^0)}$ we can show that 
$\Phi_D\circ\Psi_D$ is equivalent to the identity functor on
$\D(J)$. It follows that $\Phi_D$ is an equivalence. 
\end{proof}

\begin{corollary}
\label{singular_Ext}
Let $D$ be an arbitrary curve of type (2), (3), or (4), and let
$J:=\overline{Jac}D$ be its compactified Jacobian. Let $L_1$ and $L_2$
be two non-isomorphic torsion-free rank one sheaves of degree zero on
$J$. Then
$$\mathrm{Ext}^p_J(L_1,L_2)=0$$
for all $p$.
\end{corollary}

\begin{proof}
This follows from Corollary~\ref{singular_Mukai}, since $L_1$ and
$L_2$ are the Fourier-Mukai transforms of two skyscraper sheaves
supported at distinct points of $\overline{\mathrm{Pic}}^0J$, and
$\Phi_D$ preserves the $\mathrm{Ext}^{\bullet}$-pairing. However, we
can give a direct proof using Theorem~\ref{untwisted_equivalence},
which says that $\Phi^{\U}_{P\rightarrow Z^0}$ is an equivalence
of triangulated categories. By Bridgeland's criterion, 
Theorem~\ref{bridgeland}, it follows that we must have
$$\mathrm{Ext}^i_{Z^0}(\U_{m_1},\U_{m_2})=0$$
for all integers $i$ and all $m_1\neq m_2\in P$. In particular,
suppose that $m_1$ and $m_2$ are points in the singular fibre
$\overline{\mathrm{Pic}}^0J$ corresponding to $L_1$ and $L_2$, 
so that $\U_{m_1}=\iota_*L_1$ and $\U_{m_2}=\iota_*L_2$.

By the remark preceeding Lemma~\ref{after_flat}, the fibration
$Z^0\rightarrow\P^2$ is
flat. Therefore we can pull-back the Koszul resolution of a point in
$\P^2$ ({\em any\/} point) to get a resolution of the structure sheaf
of the corresponding fibre. It follows that the spectral sequence in 
Section 7.2 of Bridgeland and Maciocia~\cite{bm02} exists also for
singular fibres in our case, and thus we have
$$E^{p,q}_2:=\mathrm{Ext}^p_J(L_1\otimes\Lambda^q\O_J^{\oplus
  2},L_2)\Longrightarrow\mathrm{Ext}^{p+q}_{Z^0}(\iota_*L_1,\iota_*L_2).$$
The vanishing of the right hand side allows us to conclude that
$$E^{0,0}_2=\mathrm{Ext}^0_J(L_1,L_2)=\mathrm{Ext}^0_{Z^0}(\iota_*L_1,\iota_*L_2)$$
vanishes (which also follows since $L_1$ and $L_2$ are stable and
not isomorphic). Suppose 
$\mathrm{Ext}^p_J(L_1,L_2)$ vanishes for all $p\leq k$. Then
$E^{p,q}_2$ vanishes for $p\leq k$ and for all $q$. Therefore
$$E^{p+1,0}_2=\mathrm{Ext}^{p+1}_J(L_1,L_2)=\mathrm{Ext}^{p+1}_{Z^0}(\iota_*L_1,\iota_*L_2)$$
also vanishes. By induction
$$\mathrm{Ext}^p_J(L_1,L_2)$$
vanishes for all $p$. Note that since $J$ is not smooth, $L_1$ and
$L_2$ need not have finite projective resolutions. Thus we show the
vanishing for all $p$, rather than just $p\leq 2=\mathrm{dim}J$.
\end{proof}

\begin{remark}
Note that $J$ has trivial normal bundle in $Z^0$, and $Z^0$ has
trivial canonical bundle, so by adjunction the canonical bundle $\K_J$
of $J$ is also trivial. Also, the sheaves on $J$ parametrized by 
$\overline{\mathrm{Pic}}^0J$ are clearly simple. One might expect that
combining these observations with Corollary~\ref{singular_Ext}, we
could obtain another proof of Corollary~\ref{singular_Mukai}, using
Bridgeland's criterion, Theorem~\ref{bridgeland}. However, Bridgeland's
result is for smooth varieties and it is not immediately clear how to
generalize his proofs in~\cite{bridgeland99} to singular varieties. 
\end{remark}

\begin{remark}
Let us make one more remark about 
Theorem~\ref{untwisted_equivalence}: although $P\cong Z^0$, it is 
important to distinguish these spaces in the sense that $P$ should
be regarded as the dual fibration of $Z^0$. In this paper, our 
four-folds $Z^0$ and $Z^1$ are fibred by principally polarized 
abelian surfaces, which are self-dual, which is why $Z^0$ happens
to be isomorphic to its dual fibration. 

However, there exist holomorphic symplectic manifolds which are 
fibred by non-principally polarized abelian varieties. The 
generalized Kummer varieties provide examples (see 
Debarre~\cite{debarre99}); indeed Proposition~5.3 
of~\cite{sawon03} shows that generalized Kummer varieties cannot 
be fibred by principally polarized abelian varieties. The 
possibility of constructing a Fourier-Mukai transform for 
generalized Kummer four-folds was briefly discussed in Section~5.4 
of~\cite{sawon04i}; there remain many technical details to be 
resolved.
\end{remark}

\subsection{A twisted Fourier-Mukai transform} 

At the beginning of this section we made the observation that
$$P:=\overline{\mathrm{Pic}}^0(Z^0/\P^2)\cong\overline{\mathrm{Pic}}^0(Z^1/\P^2).$$ 
Thus $P$, which is by definition the dual fibration of $Z^0$, is
also the dual fibration of $Z^1$. In particular, $P$ parametrizes 
stable sheaves on $Z^1$, so there exist local universal sheaves on 
$Z^1\times P_i$ (for some cover $\{P_i\}$ of $P$). As in Section 2, 
there exists a gerbe $\beta\in\H^2(P,\O^*)$ which is the obstruction 
to the existence of a global universal sheaf on $Z^1\times P$. Since 
$Z^1$ does not admit a section, we know from 
Proposition~\ref{abelian_prop} that $\beta$ is non-zero. The 
collection of local universal sheaves gives us a 
$\pi^*_P\beta$-twisted universal sheaf $\U$ on $Z^1\times P$, where 
$\pi_P$ is projection to $P$. We can therefore construct the functor
$$\Phi^{\U}_{P\rightarrow Z^1}:\D (P,\beta^{-1})\rightarrow\D (Z^1).$$ 

\begin{theorem}
\label{twisted_equivalence}
The functor $\Phi^{\U}_{P\rightarrow Z^1}$ is an equivalence of
triangulated categories.
\end{theorem}

\begin{proof}
We will apply C{\u a}ld{\u a}raru's 
Proposition~\ref{equivalence_twisted}, which is the twisted version 
of Bridgeland's criterion, Theorem~\ref{bridgeland}. Let
$\U_m:=\Phi^{\U}_{P\rightarrow Z^1}\O_m$ be the sheaf on $Z^1$ which
the point $m\in P$ parametrizes. We must show
\begin{enumerate}
\item for all $m\in P$, $\U_m\otimes\K_{Z^1}=\U_m$ and $\U_m$ is
  simple,
\item for all integers $i$ and all $m_1\neq m_2\in P$,
$$\mathrm{Ext}^i_{Z^1}(\U_{m_1},\U_{m_2})=0.$$
\end{enumerate}
The proof of condition (1) is the same as in the proof of 
Theorem~\ref{untwisted_equivalence}. Namely, since $\U_m$ is stable, 
it is also simple. Since $Z^1$ is holomorphic symplectic, it has 
trivial canonical bundle $\K_{Z^1}$.

Regarding condition (2), first suppose that $m_1$ and $m_2$ lie in 
different fibres of $P$. Then as before $\U_{m_1}$ and $\U_{m_2}$ are 
sheaves supported on different (disjoint) fibres of $Z^1$. Therefore
$$\mathrm{Ext}^i_{Z^1}(\U_{m_1},\U_{m_2})=0$$
for all $i$ because all local ${\E}xt$ sheaves vanish.

Next suppose that $m_1$ and $m_2$ lie in the same fibre of $P$, which 
may be smooth or singular. Then $\U_{m_1}$ and $\U_{m_2}$ are of the 
form $\iota_*L_1$ and $\iota_*L_2$ respectively, where 
$\iota:Z^1_t\hookrightarrow Z^1$ is inclusion, and $L_1$ and $L_2$ 
are torsion-free rank one sheaves on $Z^1_t$ of degree one. As
before we have the spectral sequence
$$E^{p,q}_2:=\mathrm{Ext}^p_{Z^1_t}(L_1\otimes\Lambda^q\O_{Z^1_t}^{\oplus
  2},L_2)\Longrightarrow\mathrm{Ext}^{p+q}_{Z^1}(\iota_*L_1,\iota_*L_2)$$
which exists for both smooth and singular fibres (see the comments in
the proof of Corollary~\ref{singular_Ext}). If $Z^1_t$ is smooth, and
$L_1$ and $L_2$ are not isomorphic, then the spectral sequence
vanishes as in the proof of Theorem~\ref{untwisted_equivalence}. If
$Z^1_t$ is singular, then it is the compactified Jacobian of a curve
of type (2), (3), or (4). Then Corollary~\ref{singular_Ext} showed
that 
$$\mathrm{Ext}^p_{Z^1_t}(L_1,L_2)$$
vanishes for all $p$ if $L_1$ and $L_2$ are not isomorphic. So once
again, all terms in the spectral sequence vanish. It follows that
$$\mathrm{Ext}^i_{Z^1}(\U_{m_1},\U_{m_2})=0$$ 
for all $i\in\Z$ and for all $m_1\neq m_2\in P$. This concludes
the proof of condition (2), and of the theorem.
\end{proof}

\begin{remark}
Since the proof of Bridgeland and Maciocia's
Theorem~\ref{bridgeland_maciocia} only relies on local arguments, it
can be generalized to the twisted case. This would lead to a direct
proof of Theorem~\ref{twisted_equivalence}, without the need to first
prove Theorem~\ref{untwisted_equivalence} and thereby obtain
Corollary~\ref{singular_Ext}.
\end{remark}

\begin{remark}
Recall that $Z^1$ is a torsor over $Z^0\cong P$. Since $Z^0$ admits
a section and $Z^1$ does not, we could regard $Z^1$ as being a 
`twisted' version of the space $Z^0$. Thus the theorem says that the 
derived category of twisted sheaves on the `untwisted' space $Z^0$ 
is equivalent to the derived category of (untwisted) sheaves on the 
`twisted' space $Z^1$.
\end{remark}

\subsection{Deformations of fibrations}

In this final subsection we will show that $Z^1$ and $Z^0$ can be
connected by a one parameter family in their space of deformations,
which only passes through Lagrangian fibrations. This follows by
considering the subspace of $\H^2(P,\O^*)$ consisting of gerbes on
$P$ which arise from torsors over $Z^0$, and showing that it is
connected.

We start with the exponential exact sequence
$$0\rightarrow\Z\rightarrow\O\rightarrow\O^*\rightarrow 0$$
on $P$. The cohomology of each of these sheaves on $P$ can be computed
using the projection $p_P:P\rightarrow\P^2$ and the Leray spectral
sequence. Thus we have
$$E^{i,j}_2(\Z):=\H^i(\P^2,{\bf R}^jp_{P*}\Z)\Rightarrow\H^{i+j}(P,\Z)$$
and since ${\bf R}^0p_{P*}\Z_P\cong\Z_{\P^2}$ we can compute the bottom row
of the left hand side and obtain
$$\begin{array}{ccccc}
\vdots & \vdots & \vdots & & \\
{\H}^0(\P^2,{\bf R}^2p_{P*}\Z) & {\H}^1(\P^2,{\bf R}^2p_{P*}\Z) &
{\H}^2(\P^2,{\bf R}^2p_{P*}\Z) & \ldots & \\
{\H}^0(\P^2,{\bf R}^1p_{P*}\Z) & {\H}^1(\P^2,{\bf R}^1p_{P*}\Z) &
{\H}^2(\P^2,{\bf R}^1p_{P*}\Z) & \ldots & \\
{\Z} & 0 & {\Z} & 0 & \qquad{\Z}. \\
\end{array}$$
Moreover, the right-most term $\Z$ must survive the higher derivations
$d_2(\Z)$, $d_3(\Z)$, and $d_4(\Z)$, since the class generating
$\H^4(\P^2,\Z)$ can be pulled-back to give a non-trivial class in
$\H^4(P,\Z)$ (the class of a fibre).

Next we have
$$E^{i,j}_2(\O):=\H^i(\P^2,{\bf R}^jp_{P*}\O)\Rightarrow\H^{i+j}(P,\O).$$
Matsushita~\cite{matsushita00} proved that
${\bf R}^jp_{P*}\O_P\cong\Omega^j_{\P^2}$ and therefore we can compute
the left hand side precisely, obtaining
$$\begin{array}{cccc}
\vdots & \vdots & \vdots & \\
\qquad 0\qquad & \qquad 0\qquad & \qquad{\C}\qquad & \qquad\ldots\qquad  \\
0 & {\C} & 0 & \ldots \\
{\C} & 0 & 0 & \ldots \\
\end{array}$$
The spectral sequence degenerates at the $E^2$ term and gives
$$\H^k(P,\O)\cong\left\{\begin{array}{ll} {\C} & k\mbox{ even,} \\ 0 &
    k\mbox{ odd.} \\ \end{array}\right.$$

Finally we have
$$E^{i,j}_2(\O^*):=\H^i(\P^2,{\bf R}^jp_{P*}\O^*)\Rightarrow\H^{i+j}(P,\O^*).$$
In this case we know that ${\bf R}^0p_{P*}\O_P^*\cong\O^*_{\P^2}$, and
thus we can compute the bottom row of the left hand side (using the
exponential long exact sequence on $\P^2$), obtaining
$$\begin{array}{cccccc}
\vdots & \vdots & \vdots & \\
{\H}^0(\P^2,{\bf R}^2p_{P*}\O^*) & {\H}^1(\P^2,{\bf R}^2p_{P*}\O^*) &
{\H}^2(\P^2,{\bf R}^2p_{P*}\O^*) & \ldots & \\
{\H}^0(\P^2,{\bf R}^1p_{P*}\O^*) & {\H}^1(\P^2,{\bf R}^1p_{P*}\O^*) &
{\H}^2(\P^2,{\bf R}^1p_{P*}\O^*) & \ldots & \\
{\O}^* & {\Z} & 0 & {\Z} & \ldots \\
\end{array}$$

Now $Z^1$ is a torsor over $Z^0$, corresponding to a gerbe
$\beta\in\H^2(P,\O^*)$. Moreover, there is a cover $\{P_i\}$ of $P$
obtained by pulling back a cover of $\P^2$, and $\beta$ can be
represented by line bundles $\L_{ij}$ on pair-wise intersections
$P_{ij}$. Moreover, these line bundles have degree zero on each fibre
of $p_P:P\rightarrow\P^2$. Conversely, given a gerbe $\beta^{\prime}$
with these properties, we can work backwards to construct the torsor
$Z^{\beta^{\prime}}$ over $Z^0$ (see Subsection 2.4). Now a family of
line bundles on fibres of $P$ is a local section of ${\bf R}^1p_{P*}\O^*$,
and the degree is given by the coboundary map $\delta_1$ of the long
exact sequence of direct image sheaves
$$\ldots\rightarrow {\bf R}^1p_{P*}\O\rightarrow
{\bf R}^1p_{P*}\O^*\stackrel{\delta_1}{\longrightarrow}{\bf
  R}^2p_{P*}\Z\rightarrow {\bf R}^2p_{P*}\O\rightarrow\ldots$$
Thus $\beta$ corresponds to a cocycle in $\H^1(\P^2,{\bf R}^1p_{P*}\O^*)$
which moreover can be represented by local sections of
$$\mathrm{ker}\delta_1\subset {\bf R}^1p_{P*}\O^*.$$
Looking at the spectral
sequence, we see that this cocycle must also lie in the kernel of
$$d_2(\O^*):\H^1(\P^2,{\bf R}^1p_{P*}\O^*)\rightarrow\H^3(\P^2,{\bf R}^0p_{P*}\O^*)\cong\Z$$
since it survives to give a class in $\H^2(P,\O^*)$. We claim that
this makes one of the earlier conditions redundant.

\begin{lemma}
\label{degree_zero}
Suppose that $\alpha$ lies in the kernel of 
$$d_2(\O^*):\H^1(\P^2,{\bf R}^1p_{P*}\O^*)\rightarrow\H^3(\P^2,{\bf R}^0p_{P*}\O^*)\cong\Z.$$
Then $\alpha$ can be represented by line bundles $\L_{ij}$ on
pair-wise intersections $P_{ij}$ which have degree zero on each fibre
of $p_P:P\rightarrow\P^2$.
\end{lemma}

\begin{proof}
By the functoriality of spectral sequences, the coboundary maps of the
long exact sequence of direct image sheaves
$$\ldots\rightarrow {\bf R}^0p_{P*}\O^*\stackrel{\delta_0}{\longrightarrow}
{\bf R}^1p_{P*}\Z \rightarrow {\bf R}^1p_{P*}\O\rightarrow
{\bf R}^1p_{P*}\O^*\stackrel{\delta_1}{\longrightarrow}{\bf
  R}^2p_{P*}\Z\rightarrow\ldots$$
induce maps which commute with the derivations of the spectral
sequences. Thus we obtain a commutative diagram
$$\begin{array}{ccc}
{\H}^1(\P^2,{\bf R}^1p_{P*}\O^*) & \stackrel{d_2(\O^*)}{\longrightarrow} &
{\H}^3(\P^2,{\bf R}^0p_{P*}\O^*) \\
\downarrow{\H}^1(\delta_1) & & \downarrow{\H}^3(\delta_0) \\
{\H}^1(\P^2,{\bf R}^2p_{P*}\Z) & \stackrel{d_2(\Z)}{\longrightarrow} &
{\H}^3(\P^2,{\bf R}^1p_{P*}\Z). \\
\end{array}$$
Therefore since $\alpha$ is in the kernel of $d_2(\O^*)$,
$\H^1(\delta_1)\alpha$ must lie in the kernel of $d_2(\Z)$. Now as
observed above, $\H^4(\P^2,{\bf R}^0p_{P^*}\Z)\cong\Z$ must survive to
give a class in $\H^4(P,\Z)$, and thus the map
$$d_3(\Z):\mathrm{ker}d_2(\Z)\subset{\H}^1(\P^2,{\bf
  R}^2p_{P*}\Z)\rightarrow{\H}^4(\P^2,{\bf R}^0p_{P*}\Z)$$
is trivial. This means that $\H^1(\delta_1)\alpha$ also lies in the kernel
of $d_3(\Z)$, and thus it survives to give a class in
$\H^3(P,\Z)$. However, $P$ is a deformation of the Hilbert scheme of
two points on a K3 surface, and thus $\H^3(P,\Z)=0$. We conclude that
$\H^1(\delta_1)\alpha=0$.

So $\alpha$ is given by local sections $\alpha_{ij}$ of
${\bf R}^1p_{P*}\O^*$ (i.e.\ the line bundles $\L_{ij}$) over the
pair-wise intersections $P_{ij}$, such that the cocycle
$\{\delta_1\alpha_{ij}\}$ is actually a coboundary. So there exist
local sections $\gamma_i$ of ${\bf R}^2p_{P*}\Z$ over $P_i$ such that
$$\delta_1\alpha_{ij}=\gamma_i-\gamma_j.$$
Consider the image of $\gamma_i$ under the map
$${\bf R}^2p_{P*}\iota:{\bf R}^2p_{P*}\Z\rightarrow {\bf R}^2p_{P*}\O$$
induced by the inclusion $\iota:\Z\rightarrow\O$. Observe that
$${\bf R}^2p_{P*}\iota\gamma_i-{\bf R}^2p_{P*}\iota\gamma_j={\bf
  R}^2p_{P*}\iota(\delta_1\alpha_{ij})$$  
vanishes, by the exactness of
$$\ldots\rightarrow {\bf R}^1p_{P*}\O^*\stackrel{\delta_1}{\longrightarrow}
{\bf R}^2p_{P*}\Z\stackrel{{\bf R}^2p_{P*}\iota}{\longrightarrow}
{\bf R}^2p_{P*}\O\rightarrow\ldots$$
Therefore the ${\bf R}^2p_{P*}\iota\gamma_i$ agree on overlaps and can
be patched together to give a global section of ${\bf
  R}^2p_{P*}\O$. However, Matsushita proved in~\cite{matsushita00} that
${\bf R}^2p_{P*}\O\cong\Omega^2_{\P^2}$, which has no global sections,
and therefore ${\bf R}^2p_{P*}\iota\gamma_i$ vanishes for all $i$.

Using the exactness of the long exact sequence of direct image sheaves
once again, we conclude that there exist local sections $\epsilon_i$
of ${\bf R}^1p_{P*}\O^*$ over $P_i$ such that
$$\gamma_i=\delta_1\epsilon_i.$$
Now define 
$$\alpha^{\prime}_{ij}:=\alpha_{ij}-\epsilon_i+\epsilon_j.$$
These are local sections of ${\bf R}^1p_{P*}\O^*$ over $P_{ij}$ (we have
written the group action additively, but if one wants to think of
these as families of line bundles on fibres then the group action is
just tensor product). Moreover the collection
$\{\alpha^{\prime}_{ij}\}$ represents the same cohomology class in
$\H^1(\P^2,{\bf R}^1p_{P*}\O^*)$ as $\alpha$, though now we have
\begin{eqnarray*}
\delta_1\alpha^{\prime} & = &
\delta_1\alpha_{ij}-\delta_1\epsilon_i+\delta_1\epsilon_j \\
 & = & \gamma_i-\gamma_j-\gamma_i+\gamma_j \\
 & = & 0.
\end{eqnarray*}
So $\alpha^{\prime}_{ij}$ corresponds to a line bundle
$\L^{\prime}_{ij}$ on $P_{ij}$ which has degree zero on each fibre of
$p_P:P\rightarrow\P^2$. This completes the proof.
\end{proof}

Finally we prove the following result.

\begin{theorem}
\label{Lagrangian_deformation}
There is a one-parameter family of Lagrangian fibrations
$Z^t$ over $\P^2$, connecting $Z^0$ and $Z^1$ in their moduli space of
deformations. Each fibration $Z^t\rightarrow\P^2$ is a torsor over
$Z^0$, and it corresponds to a gerbe $\beta_t\in\H^2(P,\O^*)$.
\end{theorem}

\begin{proof}
Let $\beta\in\H^2(P,\O^*)$ be the gerbe corresponding to $Z^1$, and
$\alpha$ the corresponding class in $\H^1(\P^2,{\bf
  R}^1p_{P*}\O^*)$. Then $\alpha$ survives the spectral sequence to
give the class $\beta$. Part of the exponential long exact sequence on
$P$ looks like
$$\ldots\rightarrow\C\stackrel{\H^2(\mathrm{exp})}{\longrightarrow}\H^2(P,\O^*)\rightarrow\H^3(P,\Z)\rightarrow\ldots$$
since $\H^2(P,\O)\cong\C$. We know that $\beta$ can be represented
by line bundles $\L_{ij}$ on $P_{ij}$ which have degree zero on each
fibre of $p_P:P\rightarrow\P^2$. This implies that
$\H^1(\delta_1)\alpha=0$, and hence $\beta$ must map to zero in
$\H^3(P,\Z)$. Note that we have shown this without using the fact that
$\H^3(P,\Z)$ vanishes.

By exactness, $\beta$ is the image under $\H^2(\mathrm{exp})$ of an
element in $\H^2(P,\O)\cong\C$. Likewise, $\alpha$ must be the image
of a class
$$\kappa\in\H^1(\P^2,{\bf R}^1p_{P*}\O)\cong\C$$
under the map $\H^1({\bf R}^1p_{P*}\mathrm{exp})$ coming from
$\mathrm{exp}:\O\rightarrow\O^*$. Next let us define
$\alpha_t\in\H^1(\P^2,{\bf R}^1p_{P*}\O^*)$ to be the image of
$t\kappa\in\H^1(\P^2,{\bf R}^1p_{P*}\O)$ under
$\H^1({\bf R}^1p_{P*}\mathrm{exp})$. Since
$$\begin{array}{ccc}
{\H}^1(\P^2,{\bf R}^1p_{P*}\O) & \stackrel{d_2(\O)}{\longrightarrow} &
{\H}^3(\P^2,{\bf R}^0p_{P*}\O)=0 \\
\downarrow{\H}^1({\bf R}^1p_{P*}\mathrm{exp}) & &
\downarrow{\H}^3({\bf R}^0p_{P*}\mathrm{exp}) \\
{\H}^1(\P^2,{\bf R}^1p_{P*}\O^*) & \stackrel{d_2(\O^*)}{\longrightarrow} &
{\H}^3(\P^2,{\bf R}^0p_{P*}\O^*). \\
\end{array}$$
commutes, $\alpha_t$ must lie in the kernel of $d_2(\O^*)$. Thus not
only does it survive to give a gerbe $\beta_t\in\H^2(P,\O^*)$, but by
Lemma~\ref{degree_zero} it can be represented by line bundles
$\L_{ij}$ on $P_{ij}$ which have degree zero on each fibre of
$p_P:P\rightarrow\P^2$. This is precisely what is required to
construct the torsor $Z^t$ over $Z^0$, concluding the proof.
\end{proof}

\begin{remark}
Note that in the proof of Lemma~\ref{degree_zero}, the vanishing of
$\H^3(P,\Z)$ was used to show that $\H^1(\delta_1)\alpha=0$. However,
in the proof of Theorem~\ref{Lagrangian_deformation} we already know
that $\H^1(\delta)\alpha_t=0$, so the argument would work even if
$\H^3(P,\Z)$ did not vanish. So even if $\H^2(P,\O^*)$ is not
connected, the gerbes on $P$ which arise from torsors over $Z^0$ form
a connected subspace. This is significant because for the generalized
Kummer four-fold $K_4$ we have $\H^3(K_4,\Z)\cong\Z^{\oplus 8}$. In
Section 5.4 of~\cite{sawon04i} the author suggested one might be able
to find new deformation classes of holomorphic symplectic four-folds
by constructing Lagrangian fibrations from gerbes in different
connected components of $\H^2(K_4,\O^*)$. Unfortunately the above
argument implies that this will not work. 
\end{remark}

\begin{remark}
Let $\mathrm{Def}(Z^0)$ be the Kuranishi space parametrizing
deformations of $Z^0$ as a complex manifold. In~\cite{sawon04ii} the
author proved that there is a subspace
$\Delta\subset\mathrm{Def}(Z^0)$ of codimension one parametrizing
deformations of $Z^0$ which are Lagrangian fibrations, and a subspace
$\Delta^{\prime}\subset\Delta\subset\mathrm{Def}(Z^0)$ of codimension
one in $\Delta$ (and hence codimension two in $\mathrm{Def}(Z^0)$)
parametrizing deformations of $Z^0$ which are Lagrangian fibrations
with global sections. The one-parameter family described in
Theorem~\ref{Lagrangian_deformation} is consistent with these results;
it can be regarded as a deformation inside $\Delta$ but transverse to
$\Delta^{\prime}$.
\end{remark}



\begin{flushleft}
Department of Mathematics\hfill sawon@math.sunysb.edu\\
SUNY at Stony Brook\hfill www.math.sunysb.edu/$\sim$sawon\\
Stony Brook NY 11794-3651\\
USA\\
\end{flushleft}

\end{document}